\documentclass[11pt]{article}

\usepackage{latexsym}
\usepackage{amsmath,amsthm,amsfonts,amssymb}
\usepackage[T1]{fontenc}
\usepackage[latin1]{inputenc}
\usepackage{aeguill}
\usepackage{url}
\usepackage{enumerate}
\usepackage{mdwlist} 
\usepackage{rotating} 
\usepackage{pst-all} 
\usepackage{export}
\usepackage{pst-ellipse}
\usepackage[all]{xypic}

\newtheorem{thm}{Theorem}[section]

\newtheorem*{thm*}{Theorem}
\newtheorem{dfn}[thm]{Definition} 
\newtheorem*{dfn*}{Definition}

\newtheorem{cor}[thm]{Corollary}
\newtheorem*{cor*}{Corollary}

\newtheorem{prop}[thm]{Proposition} 
\newtheorem*{prop*}{Proposition} 
\newtheorem*{properties*}{Properties} 
 
\newtheorem{lem}[thm]{Lemma} 
\newtheorem*{lem*}{Lemma}

\newtheorem*{claim*}{Claim} 
 
\newtheorem*{fact*}{Fact} 
\newtheorem{fact}[thm]{Fact}

\theoremstyle{remark}
\newtheorem*{rem*}{Remark}
\newtheorem{rem}[thm]{Remark}
\newtheorem*{example*}{Example}

\newenvironment{SauveCompteurs}[1]{%
\newcommand{\monparametre}{#1}
\openexport{\monparametre_sauve}
  \Export{thm}\Export{section}\Export{subsection}\Export{subsubsection}
\closeexport}{}

\newenvironment{UtiliseCompteurs}[1]{%
\newcommand{\monparametre}{#1}
\openexport{\monparametre_aux}
  \Export{thm}\Export{section}\Export{subsection}\Export{subsubsection}
\closeexport
\PackageInfo{export}{\MessageBreak
Importations from \monparametre_sauve.xpt\MessageBreak}%
\InputIfFileExists{\monparametre_sauve.xpt}{\relax}{\relax}%
\renewcommand{\label}[1]{}
}{%
\PackageInfo{export}{\MessageBreak
Importations from \monparametre_aux.xpt\MessageBreak}%
\InputIfFileExists{\monparametre_aux.xpt}{\relax}{\relax}}

\newlength{\espaceavantspecialthm}
\newlength{\espaceapresspecialthm}
\newenvironment{specialthm}[2]{\refstepcounter{thm} 
\vskip \espaceavantspecialthm \noindent \textbf{#1\thethm #2} \itshape}%
{\normalfont \vskip \espaceapresspecialthm}

\newenvironment{specialthm*}[1]{
\vskip\espaceavantspecialthm \noindent \textbf{#1} \itshape}%
{\normalfont \vskip \espaceapresspecialthm}

\newlength{\espaceavantenonce}
\newlength{\espaceapresenonce}
\newcommand{\fontetitreun}[1]{\textbf{#1}} 
\newcommand{\fontetitredeux}[1]{\textit{#1}} 

{\normalfont \vskip \espaceapresenonce}

\newenvironment{enonce1*}[1]{
\vskip\espaceavantenonce \noindent \fontetitreun{#1} \itshape}%
{\normalfont \vskip \espaceapresenonce}

{\vskip \espaceapresenonce}

\newenvironment{enonce2*}[1]{
\vskip\espaceavantenonce \noindent \fontetitredeux{#1} }%
{\vskip \espaceapresenonce}

\newcommand{\rond}[1]{\overset{\scriptscriptstyle \circ}{#1}} 
\newcommand{\es}{\emptyset}
\renewcommand{\phi}{\varphi} 
\newcommand{\m} {^{-1}} 
\newcommand{\eps} {\varepsilon}

\newcommand {\ra} {\rightarrow}

\newcommand {\onto} {\twoheadrightarrow}
\newcommand {\xra} {\xrightarrow}

\newcommand{\actson}{\,\raisebox{1.8ex}[0pt][0pt]{\begin{turn}{-90}\ensuremath{\circlearrowright}\end{turn}}\,}

\newcommand{\ol}[1]{\overline{#1}}

\newcommand{\disjoint}{\amalg}

\newcommand{\ie} {i.~e.\ }

\newcommand {\cala} {{\mathcal {A}}}   
   
\newcommand {\calc} {{\mathcal {C}}}   
\newcommand {\cald} {{\mathcal {D}}}   
   
\newcommand {\calf} {{\mathcal {F}}}

\newcommand {\cali} {{\mathcal {I}}}

\newcommand {\calp} {{\mathcal {P}}}   
\newcommand {\calq} {{\mathcal {Q}}}   
   
\newcommand {\cals} {{\mathcal {S}}}

\newcommand{\bbR} {{\mathbb{R}}}   
\newcommand{\bbQ} {{\mathbb{Q}}}

\newcommand{\Stab} {\mathop{\mathrm{Stab}}}
\newcommand{\Fix}{\mathop{\mathrm{Fix}}}
\newcommand{\Min}{{\mathrm{Min}}}

\newcommand{\Out} {\mathop{\mathrm{Out}}}

\newcommand{\Axis} {\mathop{\mathrm{Axis}}}
\newcommand{\Char} {\mathop{\mathrm{Char}}}

\begin{document}

\title{Core and intersection number for group actions on trees.}
\author{Vincent Guirardel}
\date{\today.\\ \small Fichier \texttt{\jobname.tex}}
\maketitle

\begin{abstract}
We present the construction of some kind of {\em convex core} for the product of two actions of a group on $\bbR$-trees.
This geometric construction allows to generalize and unify the intersection number of two curves 
or of two measured foliations on a surface, Scott's intersection number of splittings, 
and the apparition of surfaces in Fujiwara-Papasoglu's construction of the JSJ splitting. 
In particular, this construction allows a topological interpretation of the intersection number analogous to the definition for curves in surfaces.
As an application of this construction, we prove that an irreducible automorphism of the free group
 whose stable and unstable trees are geometric, is actually induced a pseudo-Anosov homeomorphism on a surface.
\end{abstract}
 

Consider a surface $\Sigma$ and two homotopy classes of simple closed curves $c_1,c_2\subset\Sigma$.
Denote by $i(c_1,c_2)$ their geometric intersection number.
The nullity of the intersection number is equivalent to the possibility of isotoping the curves apart.
In terms of splittings, $i(c_1,c_2)=0$ if and only if the two splittings of $\pi_1(\Sigma)$ dual to these curves
 are \emph{compatible} (\ie have a common refinement).

In \cite{Sco_symmetry} (see also \cite{ScSw_splittings, ScSw_regular}) Scott generalized this notion of intersection number
to any pair of splittings of a finitely generated group $G$. 
This intersection number is always symmetric. 
Moreover, when edge groups of the splittings are finitely generated then this intersection number is finite,
and it vanishes if and only if the two splittings are compatible. 

By Bass-Serre theory, two splittings of a group $G$ correspond to two actions of $G$
on simplicial trees $T_1,T_2$.
In this article, we give a geometric construction of a kind of \emph{convex core} for
the diagonal action of $G$ on $T_1\times T_2$ which captures the information about
the intersection number of the corresponding splittings. 
Because of its geometric nature, this construction works naturally in the context of $\bbR$-trees.
The convexity in question here is not the CAT(0) convexity, which would give a much too large set.
The useful notion in this context is \emph{fiberwise convexity}:
a subset $E\subset T_1\times T_2$ has convex fibers if for both $i\in\{1,2\}$ and every $x\in T_i$,
$E\cap p_i\m(x)$ is convex (where $p_i:T_1\times T_2\ra T_i$ denotes the canonical projection).

\begin{specialthm*}{Main Theorem.}
Let $T_1$, $T_2$ be two minimal actions of $G$ on $\bbR$-trees having non-homothetic length functions, or being irreducible. 
Assume that $T_1$ and $T_2$ are not the refinement of a common simplicial non-trivial action.

Then there exists a subset $\calc\subset T_1\times T_2$ which is the smallest non-empty closed invariant connected subset of $T_1\times T_2$
having convex fibers. Moreover, $\calc$ is CAT(0) for the induced path-metric, and
$T_1\times T_2$ equivariantly deformation retracts to $\calc$.

We call $\calc$ the \emph{core} of $T_1\times T_2$.
\end{specialthm*}

By definition, $\calc$ is unique, and is thus invariant under automorphisms of the actions $T_1$ and $T_2$.
Moreover, this construction is symmetric by definition so $\calc(T_1\times T_2)$ is naturally isomorphic to
$\calc(T_2\times T_1)$. By contrast, the symmetry of Scott's intersection number is something that
needs a proof, and does not readily follow from the definition.

The hypotheses of the main theorem are weak, and we did not give optimal hypotheses for simplicity of the statement 
(see Proposition \ref{prop_carac_core} and Corollary \ref{cor_carac_core} for more details).  
For one-edge splittings, the assumption that $T_1$ and $T_2$ are not the refinement of a common simplicial non-trivial action
is implied by the requirement that the length functions are distinct.
Without those hypotheses, one can still give a definition for $\calc$,
but two pathologies may occur: $\calc$ may be empty, and it may fail to be connected (see sections \ref{sec_nonempty} and \ref{sec_contractibility}).
There is a remedy to the non-connectedness of $\calc$: there is a standard way to enlarge it to a connected invariant set $\Hat\calc$
with convex fibers (see section \ref{sec_augmented}). 
\\

\begin{UtiliseCompteurs}{cpt_compatible}
\begin{thm}[compare \cite{ScSw_splittings}]
   Let $T_1,T_2$ be two minimal actions of $G$ on $\bbR$-trees, such that $\calc(T_1\times T_2)\neq \es$.

Then $T_1$ and $T_2$ have a common refinement if and only if $\calc$ is $1$-dimensional.
\end{thm}
\end{UtiliseCompteurs}

The proof of this fact is very natural as $\calc(T_1\times T_2)$ itself is a common refinement of $T_1$ and $T_2$.\\

In view of the corollary above, we define the intersection number of two actions of $G$ on $\bbR$-trees 
as the covolume of $\calc$ (see definition \ref{dfn_intersection}). 
For two actions on simplicial trees with the combinatorial metric, this covolume coincides with the number of orbits
of $2$-cells in $\calc/G$, so the vanishing of the intersection number for splittings is equivalent to the
compatibility.
Of course, without any hypothesis, the action of $G$ on $\calc$ may fail to be discrete, and
in this generality, it is not clear how useful this definition of the intersection number can be.
Note however that there are intersecting cases when the actions of $G$ on $T_1$ and $T_2$ are non-discrete
whereas the action of $G$ on $T_1\times T_2$ is discrete (see the application to automorphisms of free groups in section \ref{sec_automorphisms}).
\\

Our construction of the core and intersection number generalizes and unifies several notions:
{\renewcommand{\descriptionlabel}[1]%
         {\hspace{\labelsep}\textsl{#1}}
\begin{description*}
\item[Classical and Scott's intersection number.]
Our definition of the intersection number coincides with the intersection number of two curves on a surface
and with Scott's intersection number of splittings (see example 3 in section \ref{sec_examples}, and section \ref{sec_scott}).
However, contrary to Scott's approach, we do not handle \emph{codimension one immersions }
 (\emph{almost-invariant sets} in Scott's terminology) since we need to start with actions on trees.
\item[Intersection number of measured foliations.] Given two transverse measured foliations $\calf_1$, $\calf_2$ on a surface $\Sigma$, there is a well defined intersection
number $i(\calf_1,\calf_2)$ which is the volume of the singular euclidean metric on $\Sigma$ defined by the transverse measures
of $\calf_1$ and $\calf_2$.
Our intersection number of the actions of $\pi_1(\Sigma)$ on the $\bbR$-trees $T_1,T_2$ dual to $\calf_1$, $\calf_2$ coincides with $i(\calf_1,\calf_2)$
(see example 4 in section \ref{sec_intersection}).
\item[Culler-Levitt-Shalen's core.] For two trees dual to transverse measured foliations on a surface as above, 
the core of $T_1\times T_2$  is a surface, and it is equivariantly homeomorphic to the universal cover of $\Sigma$.
In this case, Culler, Levitt and Shalen have characterized this surface as the smallest non-empty, invariant, simply-connected subset
of $T_1\times T_2$ (\cite{CLS_core}).
\item[Fujiwara-Papasoglu's enclosing groups.]
In their construction of a JSJ splitting (\cite{FuPa_JSJ}), Fujiwara and Papasoglu 
produce a surface in the product of two simplicial trees which essentially coincides with our core (Prop. \ref{prop_JSJ_FP}).
In a more general setting, their construction is not symmetric in $T_1$ and $T_2$ and produces an \emph{asymmetric core} 
(see section \ref{sec_strong_asymmetric} for a definition).

In general this asymmetric core is closely related to Scott and Swarup's (asymmetric) \emph{strong intersection numbers}:
Corollary \ref{cor_carac_strong_intersection} gives an interpretation of the strong intersection number as the number of
orbits of two cells in the asymmetric core (see Corollary \ref{cor_carac_strong_intersection}).
\end{description*}
}

A first application of our construction is a topological interpretation of the intersection number
as the minimum number of connected components of the intersection of subcomplexes representing the splitting:
\begin{UtiliseCompteurs}{cpt_interpretation}
\begin{thm}
  Assume that $Y_1,Y_2\subset X$ are two 2-sided subcomplexes, which intersect transversely, and let $T_1$, $T_2$ be the two dual trees, endowed with
the action of $\pi_1(X)$. Then $i(T_1,T_2)\leq \#\pi_0(Y_1\cap Y_2)$.

Moreover, given two non-trivial actions of a group on simplicial trees $T_1,T_2$,
there exists a complex $X$ and $Y_1,Y_2\subset X$ two 2-sided subcomplexes intersecting transversely such that 
 $i(T_1,T_2)= \#\pi_0(Y_1\cap Y_2)$.
\end{thm}
\end{UtiliseCompteurs}

The intersection number of two simple closed curves on a surface $X$ can be achieved
without changing the ambient space $X$.
In \cite[Th.~6.7]{FHS_least_area}, this result was extended to splittings dual to tori or Klein bottles in a 3-manifold
by showing that for least area representants of these submanifolds, the intersection number 
of the two induced splittings equals the number of curves in their intersection.
However, in general, one may need to change the ambient space to achieve the intersection number.

A natural question about the core is its cocompactness. In \cite{Sco_symmetry},
Scott proves that the intersection number of two splittings of a finitely generated group
over finitely generated groups have a finite intersection number. 
However, there are examples of splittings of a finitely presented group (a free group) over non-finitely generated groups
having an infinite intersection number (see Lemma \ref{lem_infinite_intersection}).
In terms of group actions on $\bbR$-trees, the finite generation of edge groups of a graph of groups means that
the corresponding Bass-Serre action is \emph{geometric}, \ie dual to a measured foliation on a finite $2$-complex (see \cite{LP}).
In this setting, we get the following finiteness result:
\begin{UtiliseCompteurs}{cpt_fund_domain}
\begin{thm}
  Let $T_1$, $T_2$ be geometric actions of a finitely generated group $G$ on $\bbR$-trees.

Then there is a set $D\subset T_1\times T_2$, which is a finite union of compact rectangles,
and such that $\calc(T_1\times T_2) \subset\ol{G.D}$.
\end{thm}
\end{UtiliseCompteurs}

This does not imply the cocompactness in general because we need to take the closure of $G.D$.
However, if $T_1$ and $T_2$ are simplicial trees, $G.D$ is automatically closed, so we get that the core 
is cocompact and that $i(T_1,T_2)$ is finite (Corollary \ref{cor_finite_intersection_simplicial}). 
We also can deduce the finiteness of the intersection number
of two geometric actions when $G$ is finitely presented (Prop. \ref{prop_finite_intersection}).\\

Finally, we give an application for automorphisms of a free group. This result is proved by showing that
the core of the product of the stable and unstable trees is almost a surface.
\begin{UtiliseCompteurs}{cpt_iwip}
\begin{cor}
  Assume that $\alpha\in\Out(F_n)$ is irreducible with irreducible powers.
Let $T_1,T_2$ be the stable and unstable actions of $F_n$ on $\bbR$-trees corresponding to $\alpha$.

If $T_1$ and $T_2$ are both geometric, then $\alpha$ is induced by
a pseudo-anosov homeomorphism of a surface with boundary.
\end{cor}
\end{UtiliseCompteurs}

The paper is organized as follows. Section \ref{sec_prelim} is devoted to basic definitions and preliminaries.
Then we give the general definition of the core $\calc$ in section \ref{sec_definition}.
In section \ref{sec_nonempty}, we study the cases where $\calc$ is empty, and we give necessary and sufficient condition
characterizing the emptiness of $\calc$.
In section \ref{sec_contractibility}, we prove that the core is contractible whenever it is connected, and
we prove that it is connected whenever $T_1$ and $T_2$ are not the refinement of a common non-trivial action on a simplicial tree.
Moreover, when $\calc$ is not connected, we show a standard way to enlarge it to get an invariant contractible set $\Hat\calc$ with
convex fibers. 
We also prove the CAT(0) property at the end of this section.
In section \ref{sec_characterization}, we prove the characterization of the core as the smallest connected non-empty closed
invariant subset with convex fibers. In section \ref{sec_compatibility}, we prove that the non-vanishing of the intersection number
is essentially the only obstruction to the compatibility of two splittings.
We give our topological interpretation of the intersection of two splittings in section \ref{sec_interpretation}.
In section \ref{sec_geometric}, we prove our finiteness result for the core of geometric actions.
We prove our application to automorphisms of free groups in section \ref{sec_automorphisms}. 
We discuss the equality of our intersection number with Scott's in section \ref{sec_scott}.
In section \ref{sec_strong_asymmetric}, we introduce an asymmetric core and relate it to Scott and Swarup's strong intersection number.
Finally in section \ref{sec_JSJ}, we relate the core with Fujiwara and Papasoglu's construction of enclosing groups.
\\

This paper was much inspired by Scott and Swarup's papers on the intersection number \cite{Sco_symmetry,ScSw_splittings}.
The construction of the core followed from an attempt of a more geometric interpretation of their definitions.

\tableofcontents

\section{Definitions and preliminaries}\label{sec_prelim}

\subsection{Basic vocabulary}
An $\bbR$-tree $T$ is a metric space in which any two points are connected by a unique topological arc,
and such that this arc is actually a geodesic.
Equivalently, an $\bbR$-tree is a geodesic metric space which does not contain any embedded topological circle.
The geodesic joining two points $a,b$ is denoted by $[a,b]$.
In an $\bbR$-tree, a subset is connected if and only if it is convex; in this case, we say that this subset is a subtree.

Consider an $\bbR$-tree $T$.
A \emph{direction} at a point $x\in T$ is a connected component of $T\setminus x$.
Note that $y,z$ are in the same direction at $x$ if and only if $[x,y]\cap[x,z]$ is not reduced to $\{x\}$.
In particular, the set of directions at $x$ corresponds to the set of germs of isometric maps from $[0,\eps]$ to $T$
sending $0$ to $x$.
A \emph{branch point} in $T$ is a point at which there are at least $3$ directions.
The metric completion $\Hat T$ of an $\bbR$-tree $T$  is still an $\bbR$-tree. 
However, one usually does not work with
complete $\bbR$-trees because it is often prefered to have minimality assumptions (see below):
Points of $\Hat T\setminus T$ are terminal points of $\Hat T$ in the sense that there is exactly one direction
based at such a point.

A ray of $T$ is an isometric embedding of $\bbR_+$ into $T$.
An \emph{end} of $T$ is an equivalence class of rays under the relation of having finite Hausdorff distance.
If $S$ is a subtree of $T$, we will denote by $\partial_\infty S\subset \partial_\infty T$ the set of ends of $S$.
\\

All actions of groups on $\bbR$-trees we consider are actions by isometries.
A group $G$ acting of $T$ is \emph{elliptic} if it fixes a point in $T$.
We say that an action of a group $G$ on an $\bbR$-tree $T$ is \emph{trivial} if every element of $G$ is elliptic.
This terminology is not really standard as the usual convention is that the action of $G$ is trivial if 
$G$ is elliptic. However, when $G$ is finitely generated, $G$ is elliptic if and only if the action of $G$ is trivial.
In general, if the action of $G$ is trivial but not elliptic, then $G$ fixes a point of the completion of $T$,
 or an end of $T$.

A group action on an $\bbR$-tree $G\actson T$ is \emph{minimal} if it has no proper invariant subtree.
When the action is non-trivial, there is a unique minimal nonempty $G$-invariant 
subtree of $T$, and this subtree is the union of translation axes of hyperbolic elements of $G$.
We denote this minimal subtree by $\Min_T(G)$.

We will denote by $\Axis_T(h)$ the axis of a hyperbolic isometry $h$ in $T$.
A \emph{positive} semi-axis of $h$ is a semi-line $A\subset Axis(h)$ such that $h.A\subset A$.
If $h$ is a hyperbolic isometry of $T$, we will denote by $\omega_T(h)$ the endpoint of $T$ defined by any \emph{positive}
semi-axis $A$ of $h$.

We denote by $l_T(g)=\min\{d(x,g.x)|x\in T\}$ the translation length of an element $g\in G$.
The action is called \emph{abelian} if there is a morphism $\phi:G\ra \bbR$ such that $l_T(g)=|\phi(g)|$.
An action is abelian if and only if there is a end of $T$ fixed by $G$.
The action is \emph{dihedral} if $T$ contains an invariant line, and some element acts as a reflection on this line.
The action is \emph{irreducible} if it not abelian, and not dihedral. Equivalently, the action is irreducible if and only
if there are two hyperbolic elements whose axes have a compact (or empty) intersection.
\\

A \emph{morphism of $\bbR$-trees} $f:T\ra T'$ is a $1$-Lipschitz map such that for
each arc $I\subset T$, there is a subdivision of $I$ into finitely many sub-intervals
on which $f$ is isometric.

A map \emph{preserving alignement} is a continuous map $f:T\ra T'$ 
such that $x\in[y,z]$ implies $f(x)\in [f(y),f(z)]$.

\begin{lem}
  Consider a continous map $f:T\ra T'$. Then following are equivalent
  \begin{enumerate*}    
  \item \label{enum_align}  $f$ preserves alignment
  \item \label{enum_convex} the preimage of every convex set is convex.
  \item \label{enum_fibre} for all $x'\in T'$, $f\m(x')$ is connected
  \end{enumerate*}
\end{lem}

\begin{proof}
  Clearly, \ref{enum_align} implies \ref{enum_convex} which implies \ref{enum_fibre}.
Now assume that \ref{enum_fibre} holds. Let $x\in [y,z]$, and assume that $f(x)\notin [f(y),f(z)]$.
Let $y'\in [y,x]$ by the point closest to $x$ such that $f(y')=f(x)$, and define $z'\in [z,x]$ similarly.
One has $y',z'\in f\m(f(x))$ and $x\notin f\m(f(x))$ contradicting the connexity of $f\m(f(x))$.
\end{proof}

The following notion of refinement generalizes the notion of refinement of a splitting to $\bbR$-trees.
\begin{dfn}[Refinement]\label{dfn_refinement}
  Consider two actions of $G$ on $\bbR$-trees $T$ and $T'$. Then one says that 
$T$ is a refinement of $T'$ if there exists an equivariant map preserving alignment
from $T$ onto $T'$.
\end{dfn}

\subsection{Technical minimality statements} \label{sec_minimality}

In what follows, given $A\subset G$, $\langle A\rangle$ denotes the group generated by $A$.

\begin{lem}\label{lem_semigroup}
Consider an action of a group $G$ on an $\bbR$-tree $T$, 
and a finitely generated semigroup $S\subset G$ acting non trivially on $T$ 
such that the minimal subtree invariant by $\langle S\rangle$ is not a line.
Let $I$ be an arc contained in the axis of a hyperbolic element $h\in S$.

Then there exists a finitely generated semigroup $S'\subset S$ such that
\begin{itemize*}
\item $\langle S'\rangle=\langle S\rangle$
\item every element $g\in S'\setminus\{1\}$ is hyperbolic in $T$, its axis contains $I$,
and $g$ translates in the same direction as $h$ on $I$.
\end{itemize*}
\end{lem}

\begin{cor}\label{cor_minimal}
  Let $T_1,T_2$ be two non-trivial actions of a group $G$ on $\bbR$-trees.
Assume that the minimal subtree of $T_1$ is dense in $T_1$.

Then for each direction $\delta_1$ in $T_1$,
there exists an element $h$ which is hyperbolic in both $T_1$ and $T_2$,
and having a positive semi-axis in $T_1$ contained in $\delta_1$.
\end{cor}

\begin{proof}[Proof of Corollary \ref{cor_minimal}]
First assume that $T_1$ contains an invariant line, which implies that $T_1$ itself is a line, and that $\delta_1$ is a semi-line.
Then $G$ has a subgroup $G_0$ of index at most $2$ consisting of elements hyperbolic in $T_1$.
If every element of $G_0$ is elliptic in $T_2$, then the action of $G$ on $T_2$ is trivial, a contradiction.
Take $h\in G_0$ which is hyperbolic in $T_2$, and $h$ or $h\m$ satisfies the conclusion of the corollary.

Now assume that $T_1$ contains no invariant line.
By density, $\delta_1$ intersects the minimal subtree of $T_1$.
Let $I_1$ be a non-degenerate arc contained in the intersection of $\delta_1$ with the translation axis of a hyperbolic element of $h\in G$.
Remember that $\omega_{T_1}(h)\in\partial_\infty T_1$ denotes the endpoint of a positive semiaxis of $h$ in $T_1$.
Up to changing $h$ to $h\m$ we can assume that $\omega_{T_1}(h)\in\partial_\infty\delta_1$.
Since $T_1$ contains no invariant line, there exists an element $h'$ whose axis is distinct from the axis of $h$.
Take $h_2,h'_2$ two elements which are hyperbolic in $T_2$ and having disting axes in $T_2$.

Apply Lemma \ref{lem_semigroup} in $T_1$ to $S=\langle h,h',h_2,h'_2\rangle$, $I$ and $h$,
to get a semigroup $S'$ such that $\langle S' \rangle=S$ and whose elements are hyperbolic in $T_1$, whose
axes contain $I$ and which translate in the same direction as $h$ on $I$.
This implies that for every $g\in S'$, $\omega_{T_1}(g)\in\partial_\infty(\delta_1)$.
Since $\langle S'\rangle=\langle S \rangle$, Serre Lemma implies that $S'$ cannot consist only of elements which are elliptic in $T_2$.
Now any $g\in S'$ which is hyperbolic in $T_2$ satisfies the conclusion of the corollary.
\end{proof}

\begin{proof}[Proof of lemma \ref{lem_semigroup}]
Let $A$ be a finite generating set of $S$ containing the hyperbolic element $h$. We apply some transformations on $A$
so that at each step, the semigroup generated by $A$ decreases, but the group generated by $A$
remains constant.


  \paragraph{Step 1: replacing $A$ by hyperbolic elements.} 
We want to replace $A$ by a set where all elements are hyperbolic.
We leave the proof of the following easy fact to the reader (the case where the intersection is empty is proved in \cite[Lem. 3.2.2]{Chi_book}):
\begin{fact}\label{fact_hyp-ell}
  Let $a,b$ be two isometries of an $\bbR$-tree $T$, with $b$ hyperbolic, and $a$ elliptic.
If $\Fix a\cap \Axis(b)$ is either empty or not reduced to a point, 
then $ab$ is hyperbolic.
\end{fact}

Thus, if $a\in A$ is elliptic, and if $\Fix a\cap \Axis(b)$ is either empty or not reduced to one point for some hyperbolic $b\in A$,
we may replace $a$ by $ab$ in $A$.

If $\Fix a\cap \Axis(b)=\{O\}$, this is more delicate. 
\begin{description}
\item[Case 1:] $a.\Axis(b)\neq \Axis(b)$.
  \begin{description}
  \item[Subcase 1.a:] $a.\omega_T(b)\neq \omega_T(b\m)$. In this case, 
one easily checks that $b^ka$ is hyperbolic for large enough $k\geq 0$, so we can replace $a$ by $b^ka$ in $A$. 
\item[Subcase 1.b:]  $a.\omega_T(b)= \omega_T(b\m)$. In this case, $a.\omega_T(b\m)\neq \omega_T(b)$
since  $a.\Axis(b)\neq \Axis(b)$. This means that the period of $\omega_T(b)$ under the action of $a$ is at least $3$
(it may be infinite), so in particular, $a^2.\omega_T(b)\neq \omega_T(b\m)$ and
 $a^3.\omega_T(b)\neq \omega_T(b\m)$.
Thus, one can replace $a$ by $\{a^2,a^3\}$ in $A$, and we are done since $a^2$ and $a^3$ either fall in subcase 1.a or
fix a semi-axis of $b$ in which case the fact \ref{fact_hyp-ell} above applies.
  \end{description}

\item[Case 2:] for every elliptic element $a\in A$ and every hyperbolic element $b\in A$ one has
$a.\Axis(b)=\Axis(b)$. 
Note that all the hyperbolic elements of $A$ cannot have the same axis since otherwise this axis would be $S$-invariant.
So consider $b,b'\in A$ having distinct axes. Then both $b^ka$ and $ab^k$ act as reflections on $\Axis(b)$ for all $k$,
and  as $k$ tends to $+\infty$, their fix points converge respectively to $\omega_T(b)$ and $\omega_T(b\m)$.
Thus, one of those fix points is not on the axis of $b'$ for some $k$. Therefore, using fact \ref{fact_hyp-ell} above, 
at least one of the two elements
$b'(b^ka)$, $b'(ab^k)$ is hyperbolic, and we can replace $a$ by this element.
\end{description}

\paragraph{Step 2: replacing $A$ by a coherent set of hyperbolic elements.}
We now aim to change $A$ so that the axis of each element  $a\in A$ contains $I$,
and so that $a$ translates in the same direction as $h$.
The lemma will follow since if $g_1,g_2$ satisfy this property, then so does $g_1g_2$.

Essentially, we are going to change $a\in A$ to some element of the form $h^kah^k$ for some large positive $k$.
The result directly follows from the following fact if no element $a\in A$ sends $\omega_T(h)$ to
 $\omega_T(h\m)$. If this occured, then neither $a^2$ nor $a^3$ would fall in this case,
and one could apply the fact after replacing $a$ by $\{a^2,a^3\}$ in $A$.
\end{proof}

\begin{fact*}
  Let $h$, $a$ be hyperbolic elements such that $a.\omega_T(h) \neq \omega_T(h\m)$.
Let $I$ be a compact interval in $\Axis(h)$.

Then for all large enough $k\geq 0$, $h^k a h^k$ is hyperbolic, its axis contains $I$,
and it translates in the same direction as $h$.
\end{fact*}

\begin{figure}[htbp]
  \centering
  \input{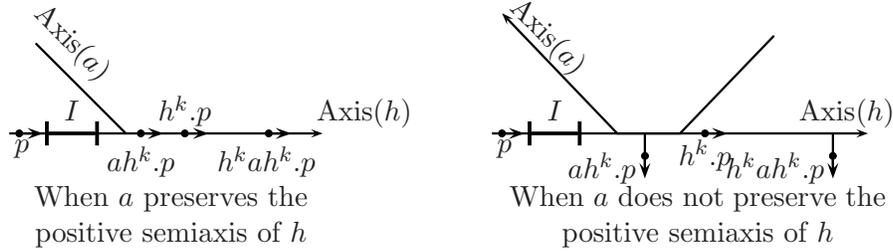}
  \caption{Making axes meet coherently}
  \label{fig_coherent}
\end{figure}

\begin{proof}[Proof of the fact.]
Of course, the fact is clear if $a$ and $h$ have the same axis.
If $a$ fixes $\omega_T(h)$, then for any given $p\in\Axis(h)$,
$h^k.p$, $ah^k.p$, and $h^kah^k.p$ are in $\Axis(h)$ for $k$ large enough,
and the fact follows easily (see figure \ref{fig_coherent}).

Otherwise, the hypothesis means that if $p$ is far enough on the positive semiaxis of $h$, then 
$a.h\notin \Axis(h)$. Figure \ref{fig_coherent} shows why the result holds in this case.
\end{proof}

\section{The main definition and examples}\label{sec_definition}

\subsection{Light quadrants and the core}

A \emph{direction} based at a point $x\in T$ is a connected component of $T\setminus x$.
A \emph{quadrant} in $T_1\times T_2$ is the product $\delta_1\times\delta_2$ of two directions $\delta_1\subset T_1$ and $\delta_2\subset T_2$.
We say that the quadrant is based at $(x_1,x_2)$ where $x_i$ is the base point of $\delta_i$.

Consider two actions of a finitely generated group $G$ 
on $\bbR$-trees $T_1$, $T_2$.
We choose a base point $*=(*_1,*_2)\in T_1\times T_2$.

\begin{dfn*}[Heavy quadrant]
 Consider a quadrant $Q=\delta_1\times\delta_2\subset T_1\times T_2$.
We say that $Q$ is \emph{heavy} if 
there exists a sequence $\gamma_k\in G$ so that 
\begin{itemize*}
\item $\gamma_k.*\in Q$
\item $d(\gamma_k.*_1,*_1)\xra{k\ra\infty}\infty$ and $d(\gamma_k.*_2,*_2)\xra{k\ra\infty}\infty$
\end{itemize*}
Otherwise, we say that $Q$ is \emph{light}.
\end{dfn*}

\begin{rem*}
This definition does not depend on the choice of the base point.
In particular, if $\delta_1$ is a direction
which does not meet the minimal subtree of $T_1$ then
for any direction $\delta_2$, $\delta_1\times \delta_2$ is light: choose the first base point in the minimal subtree
of $T_1$.
\end{rem*}

The core of $T_1\times T_2$ is what remains when one has removed the light part.
Here is a more precise definition.

\begin{dfn}
  The \emph{core} $\calc$ of $T_1\times T_2$ is the subset
$$\calc=T_1\times T_2\setminus \left[\bigcup_{Q \text{ light quadrant}} 
             Q\right]. $$
Equivalently,
$$\calc=\bigcap_{Q=\delta_1\times\delta_2 \text{ light quadrant}} 
             (\delta_1^*\times T_2 \cup T_1\times \delta_2{}^*). $$
If there is some ambiguity, we write $\calc(T_1\times T_2)$ for $\calc$.
\end{dfn}

The definition of light quadrant might seem a little bit arbitrary. 
Here are two other definitions which actually are equivalent under the weak assumption
that $\calc$ is non-empty (see Remark \ref{rem_heavy_equiv} below and Corollary \ref{cor_heavy_hyperbolic}).

\begin{dfn}[Other kinds of heavy quadrants]\label{dfn_alt_heavy}\mbox{}
  \begin{itemize*}
  \item A quadrant $Q$ is \emph{weakly heavy} if every orbit in $T_1\times T_2$ meets $Q$.
  \item A quadrant $Q=\delta_1\times \delta_2$ is \emph{made heavy by a hyperbolic element}
if there is an element $h\in G$ which is hyperbolic in $T_1$ and $T_2$, and 
such that for both $i\in\{1,2\}$, $\omega_{T_i}(h)\in\partial_\infty\delta_i$.
  \end{itemize*}
\end{dfn}

\begin{rem}\label{rem_heavy_equiv}
Clearly, a quadrant made heavy by a hyperbolic element is heavy, and a heavy quadrant is 
weakly heavy.
  If $\calc$ is non-empty, it follows conversely that any weakly heavy quadrant is heavy.
Indeed, assume that $Q$ is weakly heavy, and let $x\in\calc$. Let $g\in G$ such that
$g.x\in Q$. If $Q$ was not heavy, then $g.x$ would not lie in $\calc$, a contradiction.
\end{rem}

\subsection{Examples}\label{sec_examples}

\paragraph{Example 1: $T_1=T_2=T$, the action of $G$ on $T$ is minimal, and the set of branch points is dense.}
A quadrant  $Q=\delta_1\times \delta_2$ is light if and only if
$\delta_1\cap \delta_2=\es$. Indeed, if $\delta_1\cap\delta_2$ is empty, then $Q$ is light since
it does not meet the orbit of any point on the diagonal.
To prove the converse, assume that $\delta_1\cap\delta_2\neq\es$. 
Since branch points are dense in $T$,
the minimality of the action implies there exists a hyperbolic element $h \in G$ 
whose axis intersects $\delta_1\cap\delta_2$ in at least a semi-line, and one of the elements $h$ or $h\m$
makes $Q$ heavy.
It follows that $\calc(T\times T)$ is the diagonal of $T\times T$.

\paragraph{Example 2: $T_1=T_2=T$ is a simplicial tree with no valence $2$ vertex, and the action of $G$ on $T$ is minimal.}
The argument above extends to prove that a quadrant $\delta_1\times \delta_2$ is light if and only if
$\delta_1\cap\delta_2$ is contained in an edge.
It follows that $\calc(T\times T)=\{(v,v) | \text{$v$ vertex in $T$}\}$.
This examples illustrates the typical situation where $\calc$ is disconnected (see section \ref{sec_contractibility} for more details).

\paragraph{Example 3: $T_1$ and $T_2$ are dual to two 
non-parallel simple closed curves on a surface.}
Denote by $\Sigma$ a closed hyperbolic surface, and let $c_1,c_2$ be two distinct simple closed geodesics.
For $i\in\{1,2\}$, the tree $T_i$ dual to $c_i$ can be defined as follows:
let  $p:\Tilde\Sigma\ra \Sigma$ the universal cover of $\Sigma$; 
the vertices of $T_i$ are the connected components of $\Tilde\Sigma\setminus p\m(c_i)$
and its edges are the connected components of $p\m(c_i)$. The action of
$G=\pi_1(\Sigma)$ by deck transformations on $\Tilde \Sigma$ gives a action of $G$ on $T_i$.

We define an equivariant map $f_i:\Tilde \Sigma\ra T_i$ as follows:
Choose a small open tubular neighbourhood $A_i\simeq c_i\times ]-\eps,\eps[$ of $c_i$ foliated by curves parallel to $c_i$.
This can be done so that the two foliations of the annuli are transverse on $A_1\cap A_2$.
Let $\Tilde A_i=p\m(A_i)$, and define $f_i$ so that it sends a component of $\Tilde A_i$ to the corresponding open edge of $T_i$,
and it sends a component of $\Tilde \Sigma\setminus \Tilde A_i$ to the corresponding vertex.
Given a direction $\delta$ in $T_i$, $f_i\m(\delta)$ is at bounded Hausdorff distance from open half-plane $U_\delta$ in $\Tilde\Sigma$
bounded by a geodesic in $p\m(c_i)$.

Consider the map $F=(f_1,f_2):\Tilde\Sigma\ra T_1\times T_2$.
We are going to prove that $\calc=F(\Tilde\Sigma)$.
This will clearly follow from the fact that for any quadrant $Q=\delta_1\times \delta_2$ in $T_1\times T_2$,
$Q$ is light if and only if $f_1\m(\delta_1)\cap f_2\m(\delta_2)= \es$.
It is clear that if $f_1\m(\delta_1)\cap f_2\m(\delta_2)= \es$, then $Q$ is light since
for each point $*\in\Tilde\Sigma$, the orbit of $(f_1(*),f_2(*))$ doesn't intersect $Q$.
Conversely, if $f_1\m(\delta_1)$ intersects  $f_2\m(\delta_2)$,
then $U_{\delta_1}$ and $U_{\delta_2}$ do intersect, and there exists a element $h\in G$
whose axis in $\Tilde \Sigma$ intersects the geodesics bounding $U_{\delta_1}$ and $U_{\delta_2}$.
It is then clear that $h$ is hyperbolic in both $T_1$ and $T_2$ and that $h$ makes $Q$ heavy.

It follows that $F$ induces a bijection between the $2$-cells of $\calc$ 
and the points of $p\m(c_1)\cap p\m(c_2)$.
In other words, the number of two-cells of $\calc/G$ coincides with the intersection number $i(c_1,c_2)$. 

This observation leads to the following definition of the intersection number:

\subsection{Intersection number}\label{sec_intersection}

\begin{dfn}\label{dfn_intersection}
  Let $T_1,T_2$ be two $\bbR$-trees endowed with an action of a finitely generated
group $G$. We define the \emph{intersection number}
$i(T_1,T_2)$ as the co-volume of the action 
of $G$ on $\calc(T_1\times T_2)$ for the product measure on $T_1\times T_2$.

When $T_1$ and $T_2$ are simplicial trees with edges of length $1$, then
$i(T_1,T_2)$ is the number of $2$-cells in $\calc/G$.
\end{dfn}

Let's be more precise about this co-volume.
We say that $E\subset T_1\times T_2$ is measurable
if for every finite subtrees%
\footnote{a finite subtree is the convex hull of finitely many points.}
 $K_1\subset T_1$ and $K_2\subset T_2$,
$E\cap I_1\times I_2$ is a borel set in $K_1\times K_2$.
We denote by $\mu_{K_1,K_2}$ the product of the Lebesgue measures on $K_1$ and $K_2$.
If $E\subset T_1\times T_2$ is measurable, we define 
$$\mu(E)=\sup_{K_1,K_2} \mu_{K_1,K_2}(E\cap (K_1\times K_2)).$$
Note that a compact set may have infinite volume.

The co-volume of a measurable invariant set $\calc$ is then 
$$  \inf \Big\{\mu(E)\, | \ G.E\supset\calc,\ E\text{ measurable} \Big\}$$

\paragraph{Example 4: trees dual to two transverse measured foliations on a surface.}
Given $\calf_1,\calf_2$ two transverse measured foliations on a closed surface $\Sigma$,
one can define two $\bbR$-trees $T_i$ dual to $\calf_i$ by lifting $\calf_i$ to the universal covering
$\Tilde\Sigma$ of $\Sigma$, and taking $T_i$ the space of leaves of $\Tilde\calf_i$, endowed with the metric
given by integration of the transverse measure. 
The action of $G=\pi_1(\Sigma)$ on $\Tilde \Sigma$ provides an isometric action of $G$ on $T_1$ and $T_2$.
Denote by $f_i:\Tilde\Sigma\ra T_i$
the canonical projection, and consider $F=(f_1,f_2):\Tilde\Sigma\ra T_1\times T_2$. The argument 
of Example 3 above extends easily to this situation, so we get that $\calc(T_1\times T_2)=F(\Tilde\Sigma)$.
In particular, $\calc/G$ is isometric to $\Sigma$ endowed with the singular euclidean metric
defined by $\calf_1$ and $\calf_2$. Since the intersection number of $\calf_1$ and $\calf_2$ is the volume of
this singular metric, we get that $i(\calf_1,\calf_2)=i(T_1,T_2)$.

\subsection{Basic properties of the core}

\begin{dfn}[Convex fibers]
  Say that a subset $C\subset T_1\times T_2$ has \emph{convex fibers} (or \emph{connected fibers})
if for all $i\in\{1,2\}$ and all $x\in T_i$, $p_i\m(x)\cap C$ is convex (maybe empty).
\end{dfn}

The following properties of the core are easy:

\begin{prop}
Let $T_1,T_2$ bet two $\bbR$-trees with a non-trivial action of $G$, and let $\calc$ be the core of $T_1\times T_2$. Then
\begin{itemize*}
\item $\calc$ is closed 
\item $\calc$ has convex fibers
\item $\calc\subset \ol{\Min_{T_1}(G)}\times \ol{\Min_{T_2}(G)}$
\item if both $T_1$ and $T_2$ are simplicial, then $\calc$ is a subcomplex of $T_1\times T_2$  
\end{itemize*}
\end{prop}

\begin{rem*} If $T'_i\subset T_i$ is the minimal subtree of $T_i$, then
$\calc\subset \overline{T}'_1\times \overline T'_2$
\end{rem*}

\begin{proof}
The core is closed by definition.
  To prove that the fibers of $\calc$ are convex, just check that the fibers of the complement of a quadrant
are convex, and use the fact that an intersection of convex sets is convex.

For the third point, assume for instance that $x_1\notin \ol{\Min_{T_1}(G)}$.
Then there is a direction $\delta_1$ containing $x_1$ and not intersecting $\Min_{T_i}(G)$.
Any quadrant of the form $\delta_1\times \delta_2$ is therefore light since it does not meet 
the orbit of a base point $(*_1,*_2)$ with $*_1\in \Min_{T_1}(G)$.

Finally, assume that $T_1$ and $T_2$ are simplicial trees.
Let $Q$ be a quadrant, and let $\Hat Q$ be the union of open cells of $T_1\times T_2$ 
having a non-empty intersection with $Q$.
One easily checks that $\Hat Q$ is a quadrant,
that $T_1\times T_2 \setminus \Hat Q$ is a subcomplex of $T_1\times T_2$
(it is the product of two directions based at vertices of the trees),
and since $\Hat Q$ is contained in a bounded neighbourhood of $Q$, $\Hat Q$ is light if and only if $Q$ is light.
\end{proof}

\section{When is the core empty ?\label{sec_nonempty}}

\begin{prop} \label{prop_empty}
Let $T_1,T_2$ be two non-trivial actions of a finitely generated group $G$.
Then $\calc$ is empty if and only if 
$T_1$ and $T_2$ have homothetic length functions, 
 and
\begin{itemize*}
\item either $T_1$ (and therefore $T_2$) are dihedral 
\item or there are two ends $\omega_1$ and $\omega_2$ in $T_1$ and $T_2$ respectively,
 which are fixed by $G$, and such that 
$h$ translates towards $\omega_1$ in $T_1$ if and only if $h\m$ translates towards $\omega_2$ in $T_2$.
\end{itemize*}

In particular, if $T_1$ or $T_2$ is irreducible, then $\calc$ is non-empty.
\end{prop}

\begin{rem*}
A particular example of the second case is when $T_1$ or $T_2$ is a line, 
and $T_1$ and $T_2$ have homothetic length functions.

  If $T_1$ and $T_2$ are geometric (\ie if edge stabilizers are finitely generated in the simplicial case)
then the second case can only occur if $T_1$ and $T_2$ are both lines on which $G$ acts by translation
(see penultimate corollary in \cite{Lev_BNS}). In particular, if $T_1$ and $T_2$ are geometric and if $\calc(T_1\times T_2)=\es$,
then there is an equivariant homothety between $T_1$ and $T_2$.
\end{rem*}

We will use the two following criteria:

\begin{specialthm}{Criterion }{.}\label{crit_trois}
Assume that $a,b,c\in G$ are hyperbolic in both $T_1$ and $T_2$, and
 that for both $i\in\{1,2\}$, the three endpoints
$\omega_{T_i}(a),\omega_{T_i}(b)$ and $\omega_{T_i}(c)$ are pairwise distinct.

Let $x_i$ be the center of the triangle $\{\omega_{T_i}(a)$, $\omega_{T_i}(b)$, $\omega_{T_i}(c)\}$.
Then any quadrant containing the point $x=(x_1,x_2)$ is made heavy by a hyperbolic element.
\end{specialthm}

\begin{proof}
Consider a quadrant $\delta_1\times \delta_2$ containing $x$.
Since the complement $\delta_i^*$ of $\delta_i$ is convex, $\partial_\infty\delta_i$ must contain 
at least two of the three endpoints $\{\omega_{T_i}(a),\omega_{T_i}(b),\omega_{T_i}(c)\}$
since otherwise, $\delta_i^*$ would contain $x_i$.
Thus, there is an element $\gamma\in\{a,b,c\}$ such that for both $i\in\{1,2\}$,
$\omega_{T_i}(\gamma)\in\partial_\infty\delta_i$. Therefore, $\gamma$ makes $\delta_1\times \delta_2$ heavy.
\end{proof}

\begin{specialthm}{Criterion }{.}\label{crit_deux}
  Assume that there exist $a,b\in G$ which are hyperbolic in both $T_1$ and $T_2$, such that one of the two following hypotheses hold:
  \begin{enumerate*}
  \item either the axes of $a$ and $b$ in some $T_i$ intersect in at most one point 
  \item or the axes of $a$ and $b$ in both trees intersect in more than one point and
    \begin{itemize*}
    \item $a$ and $b$ translate in the same direction in $T_1$
    \item $a$ and $b$ translate in opposite directions in $T_2$;
    \end{itemize*}
  \end{enumerate*}

Then there is a point $x\in T_1\times T_2$ such that any quadrant containing $x$ is made heavy by a hyperbolic element.
\end{specialthm}

\begin{proof}
Consider the point $x=(x_1,x_2)$ defined as follows:
if $\Axis_{T_i}(a)\cap\Axis_{T_i}(b)\neq\es$, choose $x_i$ in this intersection, and otherwise,
take $x_i$ the projection of $\Axis_{T_i}(b)$ on $\Axis_{T_i}(a)$.
Let $E(\delta_i)=\{\gamma\in\{a,b,a\m,b\m\}|\omega_{T_i}(\gamma)\in \partial_\infty \delta_i\}$.
We are going to prove that for all quadrant $Q=\delta_1\times\delta_2$ containing $x$, $E(\delta_1)\cap E(\delta_2)\neq \es$.
It will follow that $Q$ is made heavy by an element in $E(\delta_1)\cap E(\delta_2)$, which will conclude the proof of the criterion.

For any direction $\delta_i$ containing $x_i$, 
the choice of $x_i$ implies that 
if $\Axis_{T_i}(a)\cap \Axis_{T_i}(b)$ is reduced to a point, then  $E(\delta_i)$ contains three elements.
If $\Axis_{T_i}(a)\cap \Axis_{T_i}(b)= \es$,
$E(\delta_i)$ contains either $\{a,a\m\}$, $\{a,b,b\m\}$ or $\{a\m,b,b\m\}$.
If $\Axis_{T_i}(a)\cap \Axis_{T_i}(b)\neq \es$ and if $a$ and $b$ translate in the same direction,
then  $E(\delta_i)$ contains either $\{a,b\}$ or $\{a\m,b\m\}$.
If $\Axis_{T_i}(a)\cap \Axis_{T_i}(b)\neq \es$ and if $a$ and $b$ translate in opposite directions,
then  $E(\delta_i)$ contains either $\{a,b\m\}$ or $\{a\m,b\}$.

Now the criterion follows easily:
since $\# E(T_i)\geq 2$, the only possibility allowing $E(T_1)\cap E(T_2)=\es$ is that
$E(T_1)$ and $E(T_2)$ both have two elements and are the complement of each other.
This can only occur if the axes of $a$ and $b$ have a non-degenerate intersection in $T_1$ and $T_2$.
But the hypothesis on the direction of translations prevents $E(T_1)\cap E(T_2)$ from being empty.
\end{proof}

\begin{proof}[Proof of Proposition \ref{prop_empty}.]
Let's first prove the direct implication.
We are going to exhibit a point $x\in\calc$ such that every quadrant containing $x$ is
heavy. This will imply that $x\in \calc$ so $\calc\neq\es$.

\begin{fact}\label{fact_quad1}
Assume that $T_1$ is irreducible. Then there is a point $x\in T_1\times T_2$ such that every quadrant containing
$x$ is made heavy by a hyperbolic element.
\end{fact}

\begin{proof}
  If $T_1$ is irreducible ($T_1$ is not a line and has no fix end), then one can find $a,b\in G$ which are hyperbolic in $T_1$ and whose
axes don't intersect. If $a$ and $b$ are hyperbolic in $T_2$, Criterion \ref{crit_deux} implies that $\calc\neq\es$.
Otherwise, by Corollary \ref{cor_minimal} one can find and element $h$ which is hyperbolic both in $T_1$ and $T_2$.
Now it is easy to check that there are two conjugates of $h$ by powers of $a$ or $b$ whose axes are disjoint in $T_1$:
let $p$ be the midpoint of the bridge joining $\Axis_{T_1}(a)$ to $\Axis_{T_1}(b)$, and let $\delta_a$ and $\delta_b$ be the directions
based at $p$ containg the axis of $a$ and $b$ respectively;
there are at most two integers $k$ such that $a^k.\Axis(h)\not\subset\delta_a$ (resp.\ such that $b^k.\Axis(h)\not\subset\delta_b$).
Thus Criterion \ref{crit_deux} applies to a pair $a^kha^{-k}$, $b^{k'}hb^{-k'}$.
\end{proof}

\begin{fact}\label{fact_quad2}
Assume that $T_1$ and $T_2$ are reducible (\ie are dihedral or abelian),
and that the length functions $l_1$ and $l_2$ of $T_1$ and $T_2$ are not homothetic.
  
Then there is a point $x\in T_1\times T_2$ such that every quadrant containing
$x$ is made heavy by a hyperbolic element.
\end{fact}

\begin{proof}
We restrict to a  subgroup $G'$ of $G$ of index at most $4$ so that $T_1$ and $T_2$ have a fixed end (we thus get two abelian actions).
Since a translation length satisfies $l(g^2)=2l(g)$, the length functions restricted to $G'$ are not homothetic.
The length function $l_i$ on $G'$ is the absolute value of a group morphism $\phi_i:G'\ra\bbR$.
Note that two elements $a,b$ translate in the same direction in $T_i$ if and only if $\phi_i(a)$ and $\phi_i(b)$ have the same sign.
Since $l_1$ and $l_2$ are not homothetic, consider $g,h\in G'$ such that 
$|\frac{\phi_1(g)}{\phi_1(h)}|\neq |\frac{\phi_2(g)}{\phi_2(h)}|$. 
If the signs inside the absolute values are opposite, then $g,h$ satisfy Criterion \ref{crit_deux} and we are done. 
Otherwise, up to exchanging the role of $T_1$ and $T_2$,
we can find non-zero integers $p,q$ such that
$0<\frac{\phi_1(g)}{\phi_1(h)}<p/q< \frac{\phi_2(g)}{\phi_2(h)}$.
Then $\frac{\phi_1(g^q)}{\phi_1(h^p)}<1<\frac{\phi_2(g^q)}{\phi_2(h^p)}$, 
so $\frac{\phi_1(g^q h^{-p})}{\phi_1(h^p)}<0<\frac{\phi_2(g^q h^{-p})}{\phi_2(h^p)}$,
so Criterion \ref{crit_deux} applies to $a=h^p,b=g^q h^{-p}$.
\end{proof}

\begin{fact}\label{fact_quad3}
Assume that $T_1$ and $T_2$ are reducible, have homothetic length functions, that there is no invariant line in $T_1$ nor in $T_2$,
and that for all $g\in G$, the same elements translate towards the fix end in $T_1$ and $T_2$.  

Then there is a point $x\in T_1\times T_2$ such that every quadrant containing
$x$ is made heavy by a hyperbolic element.
\end{fact}

\begin{proof}
One easily checks that one can find $g,h\in G$ whose axes in $T_1$ and $T_2$ are distinct. 
Up to changing $g$ and $h$ to their inverses, we may assume that $\omega_{T_1}(g)=\omega_{T_1}(h)$.
By hypothesis, $\omega_{T_2}(g)=\omega_{T_2}(h)$. Therefore, Criterion \ref{crit_trois} applies to  $a=g,b=g\m,c=h\m$.
\end{proof}

This ends the proof of the direct implication of the Proposition.
We now prove that $\calc$ is empty in the exceptional cases.
First, if $T_1$ is dihedral and if $T_1$ and $T_2$ have homothetic length functions,
then $T_2$ is also dihedral, and one easily checks that $\calc$ is empty.
Therefore, let $\omega_1,\omega_2$ be two ends in $T_1$ and $T_2$ respectively, which are fixed by
$G$, and such that 
$h$ translates towards $\omega_1$ in $T_1$ if and only if $h\m$ translates towards $\omega_2$ in $T_2$.
In other words, $\omega_{T_1}(h)=\omega_1$  if and only if $\omega_{T_2}(h\m)=\omega_2$. 
We prove that for any quadrant $Q=\delta_1\times\delta_2$ such that $\omega_i\notin \partial_\infty\delta_i$, then
$Q$ is light. This will clearly imply that $\calc=\es$.
Let $x_i$ be the point at which $\delta_i$ is based, and take $x=(x_1,x_2)$ as a basepoint.
Consider $g\in G$ such that $g.x_1\in \delta_1$. Since $\omega_1\in\partial_\infty\Char(g)$,
one gets that $x_1\in\Char(g)$: otherwise, the subtree $\ol{\delta_1}$
would not intersect $\Char(g)$ and would thus be disjoint from its image 
under $g$. In particular, $g$ is hyperbolic translates away from $\omega_1$.
The symmetric argument in $T_2$ says that if $g.x_2\in\delta_2$, then $g$ translates away $\omega_2$.
Therefore, $\delta_1\times\delta_2$ does not meet the orbit of the basepoint, and $\delta_1\times\delta_2$ is light.
Thus $\calc$ is empty, and this ends the proof of the proposition.
\end{proof}

\begin{rem}\label{rem_coeurvide}
Assume that $T_1$ and $T_2$ are simplicial, and given an oriented edge $e\in T_i$
denote by $\delta(e)$ the direction based at the origin of $e$ and containing the terminus of $e$.
Then the end of the proof shows that if $T_1$ and $T_2$ satisfy the second hypothesis of the proposition,
then for each pair of non-oriented edges $e_1,e_2$, there is a choice of orientations of $e_1,e_2$ (namely,
not pointing towards $\omega_i$)
and choice of base point $*$ (namely $*_i$ is the origin of $e_i$), 
such that the orbit of $*$ does not meet $\delta(e_1)\times\delta(e_2)$.
This fact implies that Scott's intersection number is also zero in this case (see section \ref{sec_scott}).
\end{rem}

\begin{cor}\label{cor_heavy_hyperbolic}
Let $T_1,T_2$ be two actions of a finitely generated group $G$ on $\bbR$-trees. Assume that $\calc(T_1\times T_2)\neq \es$.

Then each heavy quadrant is made heavy by a hyperbolic element.
\end{cor}

\begin{proof}
It follows from one of the three facts above that 
there is a point $x\in \calc$ such that 
each quadrant containing $x$ is made heavy by a hyperbolic element.
Now consider any heavy quadrant $Q=\delta_1\times\delta_2$. Since $Q$ is heavy, there is a translate $g.x$ of $x$ contained in $Q$.
Since $x\in g\m.Q$, consider a hyperbolic element $h$ making $g\m.Q$ heavy.
Then  $ghg\m$ makes $Q$ heavy.
\end{proof}

\section{Contractibility of the core}\label{sec_contractibility}

In all this section, we assume without loss of generality that $T_1$ and $T_2$ have a dense minimal subtree.
The goal of this section is to understand when the core can fail to be connected.
We saw that this occured when $T_1=T_2$ is a simplicial tree without valence 2 vertex
since in this case, $\calc=\{(v,v)|v \text{ vertex in }T\}$.
We will see that this is essentially the only case when it happens: Proposition \ref{prop_carac_connectedness}
claims that this happens if and only if $T_1$ and $T_2$ refine a common simplicial $G$-tree $T'$
as in the following diagram:
$$
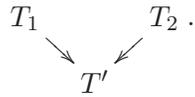
\xymatrix@=2ex{T_1\ar[dr]&&T_2\ar[dl]\\
&T'&}.$$

We will then prove that except in this pathological case, $\calc$ is contractible and that there is an equivariant
retraction by deformation of $T_1\times T_2$ onto $\calc$ (Proposition \ref{prop_contractible}).
In general, we will prove that one can enlarge $\calc$ into a natural $G$-invariant subset $\Hat\calc$ (the \emph{augmented core})
by adding some \emph{diagonals}, so that $\Hat\calc$ is contractible.

\subsection{Twice-light rectangles}\label{sec_twice-light}

\begin{figure}[htbp]
  \begin{center}
    \input{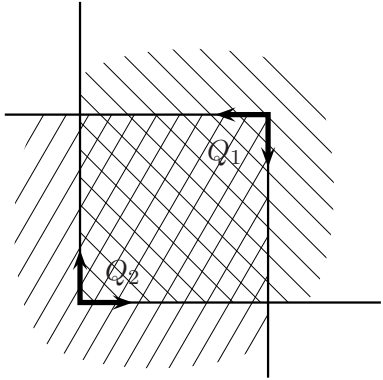}
    \caption{Two facing quadrants intersect in a twice-light rectangle}
    \label{fig_twicelight} \end{center}
\end{figure}

We are interested in the case where a point in $T_1\times T_2$
is removed \emph{twice} from $\calc$ like in figure \ref{fig_twicelight}.
This phenomenon occured in the example where $T_1=T_2$ is a simplicial tree.

We need a little bit of terminology.

\begin{dfn}[Facing directions and facing quadrants]
Say that two directions $\delta,\delta'\subset T$, based at two distinct points $x,x'$
face each other if one of the following equivalent conditions hold:
\begin{enumerate*}
\item $\delta\cup \delta'=T$
\item $\delta^*\cap \delta'{}^*=\es$
\item $]x,x'[\subset \delta\cap\delta'$
\end{enumerate*}

Say that the quadrants $Q=\delta_1\times\delta_2\subset T_1\times T_2$ and $Q'=\delta'_1\times\delta'_2\subset T_1\times T_2$ 
are \emph{facing each other} if for both $i\in\{1,2\}$, $\delta_i$ faces $\delta'_i$.
\end{dfn}

The proof of the statements contained in this definition is straightforward and left to the reader.

\begin{dfn}[Twice-light points and rectangles]
%
Consider $Q,Q'$ two light quadrants facing each other.
Then we say that $Q\cap Q'$ is a \emph{twice-light rectangle}.
\end{dfn}

Consider $Q,Q'$ two light quadrants facing each other, and let $x=(x_1,x_2)$ and $x'=(x'_1,x'_2)$
their base points. Let $R$ be the rectangle $]x_1,x'_1[\times]x_2,x'_2[$. 
The fact that $Q$ and $Q'$ face each other means that $R\subset Q\cap Q'$.
We are going to prove that actually, $R=Q\cap Q'$.

\begin{prop}[A twice light rectangle is a rectangle]\label{prop_twice-light}
Assume that the minimal subtrees of $T_1$ and $T_2$ are dense in $T_1$ and $T_2$.
Consider $Q,Q'$ two light quadrants facing each other, and let $x=(x_1,x_2)$ and $x'=(x'_1,x'_2)$
their base points. 
Let $Q\cap Q'$ be the corresponding twice light rectangle, and let $R$ be the rectangle $]x_1,x'_1[\times]x_2,x'_2[$.

Then $R=Q\cap Q'$, and none of the intervals $]x_i,x'_i[=p_i(R)$ contains a branch point.
\end{prop}

\begin{proof}
Since $Q$ and $Q'$ face each other, $R\subset Q\cap Q'$.
Now consider $(t_1,t_2)\in R$, assume that there at least three directions at $t_1$
and argue towards a contradiction (a symmetric argument will apply to $t_2$).
Let $\eta_1$ be a direction at $t_1$ which does not contain the base points of $\delta_1$ and $\delta'_1$.
Thus, $\eta_1 \subset \delta_1\cap\delta'_1$.
Since the minimal subtree of $T_1$ is dense in $T_1$, Corollary \ref{cor_minimal} 
implies that there exists an element $h\in G$
which is hyperbolic in $T_1$ and $T_2$ and 
such that $\omega_{T_1}(h)\in\partial_\infty\delta_1$.
In particular, for $k$ large enough, 
$h^k.*_1\in \eta_1$. 
Since $\delta_2\cup\delta'_2=T_2$, up to taking a subsequence, 
we may assume for instance that for $k$ large enough, $h^k.*_2\in \delta_2$.
Since $\eta_1\subset\delta_1$, this contradicts the fact that $\delta_1\times\delta_2$ is light.
\end{proof}

\subsection{Connectedness without twice light rectangles}

 Let's start with the following remark:
if $\calc$ has some twice light rectangles, then $\calc$ cannot be connected (at least if it is non-empty).
Indeed, let $Q,Q'$ be two light quadrants facing each other based at $x=(x_1,x_2)$ and $x'=(x'_1,x'_2)$ respectively.
Then $Q\cup Q'$ contains $]x_1,x'_1[\times T_2$, so $p_1(\calc)$ does not meet $]x_1,x'_1[$.
However, if $\calc$ was connected, then $p_1(\calc)$ would be a non-empty connected $G$-invariant subset of $T_1$, 
so it would be dense by our minimality hypothesis, a contradiction.

We now prove that this is the only obstruction to the connectedness of $\calc$.

\begin{prop}\label{prop_connectedness_simple}
Let $T_1,T_2$ be two actions of $G$ on $\bbR$-trees, with dense minimal subtrees, 
such that $\calc(T_1\times T_2)\neq \es$.

Then $\calc$ is connected if and only if $\calc$ 
 has no twice light rectangle. 

In this case, for each rectangle $R$ of $T_1\times T_2$,
the trace of $\calc$ on $R$ is either empty or connected.
\end{prop}

\begin{proof}
We already proved that the presence of twice light squares prevents the connectedness of $\calc$.

  Let $x,x'\in R\cap\calc$, and let $R_0$ be the rectangle $[x_1,x'_1]\times[x_2,x'_2]\subset R$.
We denote by $a=(x_1,x'_2)$ and $b=(x'_1,x_2)$ the two other corners of $R_0$.
We are going to prove that $R_0\cap \calc$ is connected, which will prove the proposition.
Consider a light quadrant $Q$ which meets $R_0$.
Since $Q$ cannot contain $x$ nor $x'$, its trace on $R_0$ is an open rectangle containing $a$ or $b$.
Moreover, if $Q$ and $Q'$ are light quadrants which intersect and such that  $a\in Q$ and $b\in Q'$, then $Q$ and $Q'$ face each other.
Since there is no twice light rectangle, this means that any light quadrant containing $a$ does not intersect a light quadrant
containing $b$.

Let $A\subset R_0$ (resp.\ $B$) be the union of the traces on $R_0$ of the light quadrants containing $a$ (resp.\ $b$).
Write $A$ (resp.\ $B$) as an increasing union $A=\cup_k A_k$ where $A_k$ is a finite union of traces of light quadrants
containing $a$. Since $\calc$ has no twice light rectangle, then for each $k$, $A_k\cap B_k=\es$.
Clearly, $R_0\setminus A_k$ and $R_0\setminus B_k$ are contractible (they are star-shaped), 
 so $R_0\setminus A_k\cup B_k$ is contractible.
Thus $R_0\setminus A\cup B$ is connected as a decreasing intersection 
of compact connected spaces.
\end{proof}

The proposition can be reformulated in a more general setting, which will be useful later:

\begin{dfn}[Coherent family of quadrants, Core]
 A family of quadrants $\calq$ of $T_1\times T_2$ is called \emph{coherent} 
if it contains no pair of quadrants facing each other.

Its \emph{core} $\calc_\calq$ is defined by $\calc_\calq=T_1\times T_2 \setminus \cup_{Q\in\calq} Q$.
\end{dfn}

\begin{prop}\label{prop_connexity_coherent}
  The core $\calc_\calq$ of a coherent family of quadrants $\calq$ is connected (maybe empty). 
Moreover, for each rectangle $R$ of $T_1\times T_2$,
the trace of $\calc_\calq$ on $R$ is either empty or connected.
\end{prop}

\subsection{The corners of a twice light rectangle}

We now study in more detail twice light rectangles in order to define the augmented core $\Hat\calc$.

\begin{lem}
Let $T_1,T_2$ be two actions of $G$ on $\bbR$-trees, such that $\calc(T_1\times T_2)\neq \es$.
Then any twice light rectangle is contained in a unique maximal twice light rectangle.
\end{lem}

\begin{proof}
  We first rule out the case where either $T_1$ or $T_2$ is a line, since in this case,
the existence of a twice light rectangle implies that $\calc(T_1\times T_2)=\es$.

Let $R=]x_1,x'_1[\times ]x_2,x'_2[$ be a twice light rectangle.
Let $]y_i,y'_i[\subset T_i$ be the maximal open interval containing $]x_i,x'_i[$, and containing no branch point 
(equivalently, $]y_i,y'_i[$ is open in $T_i$ and $y_i,y'_i$ are branch points).
Then it is clear that $]y_1,y'_1[\times]y_2,y'_2[$ is twice light, and that it is maximal for this property
because of proposition \ref{prop_twice-light}.
\end{proof}

\begin{dfn}
  Two quadrants $Q=\delta_1\times\delta_2$ and $Q'=\delta'_1\times\delta'_2$ weakly face each other 
if $\delta_1$ faces $\delta'_1$ or $\delta_2$ faces $\delta'_2$.
\end{dfn}

\begin{lem}\label{lem_facing}
Let $T_1,T_2$ be two actions of $G$ on $\bbR$-trees, whose minimal subtrees are dense.
  Assume that $Q$ and $Q'$ are two light quadrants having a non-empty intersection.
If $Q$ weakly faces $Q'$, then $Q$ faces $Q'$.
\end{lem}

\begin{proof}
Let $Q=\delta_1\times\delta_2$ and $Q'=\delta'_1\times\delta'_2$, and assume that $\delta_1$ faces $\delta'_1$.
We need to prove that $\delta_2$ faces $\delta'_2$. If this wasn't the case, then $\delta_2$ and $\delta'_2$
would be nested since $\delta_2\cap\delta'_2\neq\es$. So we assume for instance that $\delta_2\subset\delta'_2$.
Take a an element $\gamma\in G$ which is hyperbolic in $T_1$ and $T_2$ and such that $\omega_{T_2}(\gamma)\in\partial_\infty\delta_2$
(Corollary \ref{cor_minimal}). Since $\delta_1\times \delta_2$ is light, 
$\gamma^k.*_1\notin \delta_1$ for $k$ large enough, and since $\delta'_1\times \delta'_2$ is light,
$\gamma^k.*_1\notin \delta'_1$ for $k$ large enough. 
But since $\delta_1$ faces $\delta'_1$, $\delta_1\cup\delta'_1=T_1$, a contradiction.
\end{proof}

\begin{lem}
Let $T_1,T_2$ be two actions of $G$ on $\bbR$-trees, whose minimal subtrees are dense, and such that $\calc(T_1\times T_2)\neq \es$.
Let $R$ be the closure a maximal twice light rectangle. Then $\calc\cap R$ consists of exactly two opposite corners of $R$.
\end{lem}

\begin{proof}
  Let $Q=\delta_1\times\delta_2$ and $Q'=\delta'_1\times\delta'_2$ be two light quadrants facing each other 
such that $Q\cap Q'$ is the interior of $R$.
Let $(x_1,x_2)$ and $(x'_1,x'_2)$ the base points of $Q$ and $Q'$.
We prove that $R\cap \calc$ consists of the two points $a=(x_1,x'_2)$ and $b=(x'_1,x_2)$.
Clearly, $R\cap \calc\subset \{a,b\}$ since $R\setminus (Q\cup Q')=\{a,b\}$.
Now assume that $a\notin\calc$, and let $P=\eta_1\times\eta_2$ be a light quadrant and
containing $a$.

Since $R$ is the closure of a maximal twice light rectangle, both coordinates of $a$ and $b$ are
branch points. This implies that $P$ cannot face $Q$ since otherwise,
$P\cap Q$ would be a twice light rectangle, but $p_1(P\cap Q)$ would contain the branch point $p_1(a)$,
a contradiction. Similarly, $P$ cannot face $Q'$.

On the other hand, $P$ weakly faces $Q$ or $Q'$; more precisely,
one has that $\eta_1$ faces $\delta_1$
or $\delta'_1$. Indeed, $\delta_1^*\cap\delta'_1{}^*=\es$ since $\delta_1$ faces $\delta'_1$.
Now if $\eta_1$ does not face $\delta_1$ nor $\delta'_1$, then $\eta_1^*$ intersects the two subtrees
$\delta_1^*$ and $\delta'_1{}^*$, and must therefore contain the bridge joining them,
which contradicts the fact that $P$ contains $a$.

Since $P$ contains a neighbourhood of $a$, $P$ intersects both $Q$ and $Q'$.
Therefore, by lemma \ref{lem_facing}, we deduce that
$P$ faces $Q$ or $Q'$, a contradiction.
\end{proof}

\subsection{The augmented core}\label{sec_augmented}

In this section, we define the augmented core, and prove its connectedness by showing that it is the core
of a coherent family of quadrants.

\begin{dfn}
Let $T_1,T_2$ be two actions of $G$ on $\bbR$-trees, such that $\calc(T_1\times T_2)\neq \es$,
and let $R$ be a maximal twice light rectangle. The \emph{main diagonal} $\Delta_R$ is the diagonal of $R$
joining its two corners lying in $\calc$.

The \emph{augmented core} $\Hat\calc$ of $T_1\times T_2$ is the union of $\calc$ and of the main diagonal
of its maximal twice light rectangles.
\end{dfn}

\begin{prop}\label{prop_connectedness_augmented}
Let $T_1,T_2$ be two actions of $G$ on $\bbR$-trees, such that $\calc(T_1\times T_2)\neq \es$.
  Then $\Hat\calc$ is the core of a coherent family of quadrants. In particular, $\Hat\calc$ is connected,
and it intersects each rectangle into a connected set.
\end{prop}

\begin{proof}  
Let $\calq$ be the family of quadrants of $T_1\times T_2$ which don't intersect $\Hat\calc$.
Note that each quadrant of $\calq$ is light as it does not intersect $\calc$.
Moreover, $\calq$ is coherent. Indeed, assume that $Q,Q'\in\calq$ face each other. 
Then $Q\cap Q'$ is a twice light rectangle, and let $R$ be the maximal twice light rectangle containing $Q\cap Q'$. 
Since the trace of $Q\cup Q'$ on $R$ separates its two main corners,
one deduces that  $Q\cup Q'$ intersects the main diagonal $\Delta_R$, a contradiction.

There remains to prove that $\Hat\calc$ is the core of $\calq$.
We will prove the following fact later:
\begin{fact}\label{fact_un_seul_twicelight}
  A light quadrant $Q$ meets at most one maximal twice light rectangle $R$. Moreover, in this case,
the basepoint of $Q$ lies in $\ol R$.
\end{fact}

By definition, $\Hat\calc\subset \calc_\calq = T_1\times T_2 \setminus \bigcup_{Q\in\calq} Q$.
Moreover, it is clear that if $R=Q\cap Q'$ is a maximal twice light rectangle, then any element $x\in \ol{R}\setminus \Delta_R$
lies in a quadrant contained in $Q$ or $Q'$ and which does not meet $\Delta_R$.

\begin{figure}[htbp]
  \centering
  \input{quadrants_of_P.pst} 
  \caption{The quadrants of $\calp$.}
  \label{fig_quadrants_of_P}
\end{figure}

Now, let $x\notin \Hat\calc$ and which doesn't lie in the closure of a twice light rectangle. 
Let $Q=\delta_1\times\delta_2$ be a light quadrant containing $x$.
We will prove that $x$ lies in a light quadrant which does not meet any twice light rectangle.
Assume that $Q$ intersects a maximal twice light rectangle $R$.
The fact above claims that the base point of $Q$ lies in $\ol R$.
Since $R$ contains no branch point,
$Q\setminus \ol R$ can be written as the following union of quadrants
$$Q\setminus \ol R=\bigcup_{P\in \calp} P$$
where $\calp$ is defined as follows (see figure \ref{fig_quadrants_of_P}):
let $D_i$ be the set of connected components of $\delta_i\setminus p_i(\ol R)$ 
(these connected components are actually directions), and let
$\calp$ be the set of quadrants defined by
$$\calp=\{\delta_1\times \eta_2, \eta_1\times \delta_2\,|\, \eta_1\in D_1,\eta_2\in D_2\}.$$
Since $x\in P$ for some $P\in  \calp$, 
there remains to check that the quadrants of $\calp$ are light and don't intersect a twice light rectangle.

The quadrants of $\calp$ are clearly light since they are contained in $Q$.
Moreover, if a twice light rectangle $R'$ intersects a quadrant $P\in \calp$, then it also intersects $Q$,
so $R=R'$ a contradiction. Therefore, the quadrants of $\calp$ don't meet any twice light rectangle, and hence don't intersect $\Hat\calc$.
In other words, $\calp\subset\calq$ and $x\notin \calc_Q$.
\end{proof}

\begin{proof}[Proof of Fact \ref{fact_un_seul_twicelight}]
  Let $R$ be a maximal twice light rectangle contained in a light quadrant $Q=\delta_1\times\delta_2$.
We only need to prove that the base point $b=(b_1,b_2)$ of $Q$ lies in $\ol R$.
Indeed, it follows that two maximal twice light rectangles contained in $Q$ have a nonempty intersection
(they have the same germ at $b$ since they intersect the same quadrant) and must therefore coincide.

We have to prove that for both $i\in\{1,2\}$, $b_i\in p_i(\ol R)$.
So assume for instance that, $b_1\notin p_1(\ol R)$.
Then by connexity, either $p_1(\ol R)\subset \delta_1$ or  $p_1(\ol R)\subset \delta_1^*$.
The latter being impossible since $R$ meets $Q$, it follows that $p_1(R)\subset \delta_1$.
Now since $R$ meets $Q$, $\delta_2$ contains at least one of the endpoints $p_2(Q)$, call it $a_2$.
The segment $p_1(R)\times\{a_2\}$ is contained in $Q$,
and contains one of the endpoints of the main diagonal, which lies in $\calc$.
This contradicts the fact that $Q$ is light.
\end{proof}

\subsection{Characterization of the connectedness of the core}

We now can prove that the core is disconnected if and only if $T_1$ and $T_2$ 
both refine a common splitting:

\begin{prop}[Characterization of the connectedness of the core.] \label{prop_carac_connectedness}
  Let $T_1$, $T_2$ be two actions of $G$ on $\bbR$-trees, whose minimal subtrees are dense, and such that $\calc(T_1\times T_2)\neq \es$.

Then $\calc(T_1\times T_2)$ is disconnected if and only if
$T_1$ and $T_2$ are refinements of a 
common non-trivial action on a simplicial tree $T'$
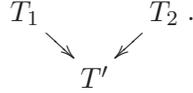
$$\xymatrix@=2ex{T_1\ar[dr]&&T_2\ar[dl]\\
&T'&}.$$
\end{prop}

\begin{rem}
The connectedness of $\calc$ is also equivalent to the absence of twice light squares (Proposition \ref{prop_connectedness_simple}),
which is equivalent to the equality $\calc=\Hat\calc$.
\end{rem}

We first prove the following lemma.
\begin{lem}\label{lem_un_par_fibre}
Given any $x_1\in T_1$, there is at most one maximal twice light rectangle $R$ such that
$x_1\in p_1(R)$.

In particular, $p_{1|\Hat\calc\setminus\calc}$ is one-to-one.
\end{lem}

\begin{proof}
  Let $R=]a_1,b_1[\times]a_2,b_2[$ and  $R'=]a'_1,b'_1[\times]a'_2,b'_2[$ be two maximal twice light rectangles as in the lemma.
We know that $a_1,b_1,a'_1,b'_1$ are branch points, and that $]a_1,b_1[$ and $]a'_1,b'_1[$ don't contain
any branch point. Since  $]a_1,b_1[$ and $]a'_1,b'_1[$ have nonempty intersection, it follows that
$]a_1,b_1[=]a'_1,b'_1[$. 
By a similar argument, if $]a_2,b_2[$ intersects $]a'_2,b'_2[$ , then $R=R'$ and we are done.
Otherwise, $]a'_2,b'_2[$ is contained in $p_2(P)$ or in $p_2(Q)$, so $R'$ is contained in $P$ or $Q$.
By lemma \ref{fact_un_seul_twicelight}, $R=R'$. 
\end{proof}

\begin{proof}[Proof of Proposition \ref{prop_carac_connectedness}]
If $T_1$ and $T_2$ both refine an action of $G$ on a simplicial tree $T'$, 
an argument similar to example 1 in section \ref{sec_examples}
shows to the existence of twice light rectangles as follows.
First, one can first assume without loss of generality that $T_1$ and $T_2$ are minimal.
Let $e=]a,a'[$ be an open edge of $T'$, and let $I_i=f_i\m(e)$.
We prove that $I_i$ contains no branch point so $I_i$ is an open interval.
Otherwise, let $b_i\in I_i$ be a branch point, and let $\delta_i$ be a direction based at $b_i$ which does not intersect $f_i\m(a)$ and $f_i\m(a')$;
first, $f_i(\delta_i)=\{f_i(b_i)\}$ since otherwise, the preimage of some point of $e$ would be disconnected.
Now, by minimality, there is an element $h\in G$ which is hyperbolic in $T_i$ and $\omega_{T_i}(h)\in\partial_\infty(\delta)$.
Now the image of every point under sufficiently high powers of $h$ ends in $\delta$.
This implies that the whole tree $T_i$ is mapped to $\{f_i(b_i)\}$ under $f_i$,
contradicting the non-triviality of the action of $G$ on $T'$.

Denote by $a_i$ and $a'_i$ the enpoints of $\ol{I_i}$ which are mapped to $a$ and $a'$ respectively under $f_i$,
and let $\delta_i$ and $\delta'_i$ be the directions in $T_i$ based at $a_i$ and $a'_i$, and containing
$I_i$. Let $Q=\delta_1\times\delta'_2$ and $Q'=\delta'_1\times\delta_2$. Choose a base point $*=(*_1,*_2)$ such that
$f_1(*_1)=f_2(*_2)\in e$. The quadrant $Q$ is light because if
$g.*\in Q$, then $g.f_1(*_1)=g.f_2(*_2)\in f_1(\delta_1)\cap f_2(\delta'_2)=e$
so $g.*\in I_1\times I_2$ and cannot go to infinity.
Since the same argument applies to $Q'$, this concludes to the existence of a twice light rectangle in $T_1\times T_2$.

To prove the converse, the idea is to obtain $T'$ by collapsing everything but the main diagonals of twice light rectangles.
More precisely, let $T'$ be the simplicial tree defined as follows: its vertices are the connected components of $\calc$,
its edges are the main diagonals of the twice light rectangles. The endpoints of edges are the natural ones.

Lemma \ref{lem_un_par_fibre} implies that each edge disconnects $T'$, so $T'$ is a tree.
Now let $x\in T_i$, and let's define $f_i(x)$. If $x\in p_i(\calc)$, then $p_i\m(x)$ is connected, so it defines a point
of $T'$ which we assign to $f_i(x)$. If $x\in p_i(R)$ for a maximal twice light rectangle $R$,
we send $x$ to the corresponding point of $\Delta_R$.
The map $f_i$ is now defined on $p_i(\calc)$ which is a dense subtree of $T_i$.
For the combinatorial metric on $T'$, $f_i$ is Lipschitz and has therefore a unique Lipschitz extension to $T_i$.
It is clear that $f_i$ preserves alignment: if $x$ separates $y$ from $z$ in $T_i$ and $f_i(x)$ is distinct
from $f_i(y)$, $f_i(z)$, then $f_i(x)$ separates $f_i(y)$ from $f_i(z)$ in $T'$.
\end{proof}

\subsection{Flow and contractibility of the augmented core}

\begin{prop}\label{prop_contractible}
  Let $T_1,T_2$ be a pair of minimal actions of $G$ on $\bbR$-trees whose core is non-empty.
Let $\Hat\calc$ be the augmented core.

Then there is an equivariant retraction of $T_1\times T_2$ onto $\Hat\calc$.
In particular, $\Hat\calc$ is contractible.
Moreover, for any rectangle $R\subset T_1\times T_2$, $R\cap\Hat\calc$ is empty or contractible.
\end{prop}

\begin{figure}[htbp]
\centering
\input{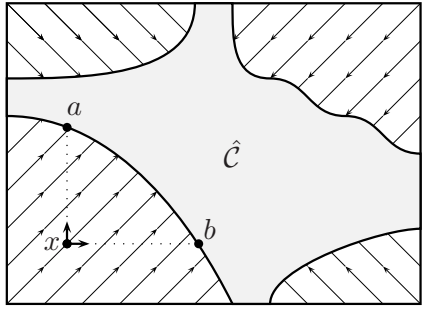}
    \caption{The equivariant flow $\phi_t$}
    \label{fig_vectorfield}
\end{figure}

\begin{rem*}
  Minimality can be replaced here by the weaker assumption that $T_i=p_i(\Hat\calc)$.
 In particular, the contractibility of $\Hat\calc$ remains true without minimality hypothesis:
$\Hat\calc$ is a deformation retract of $p_1(\Hat\calc)\times p_2(\Hat\calc)$.
\end{rem*}

We will need a slightly more general statement.

\begin{lem}\label{lem_contractile_general}
  Let $\calq$ be a family of quadrants whose core $\calc_\calq$ is connected and such that
for both $i\in\{1,2\}$, $p_i(\calc_\calq)=T_i$.

Then there is a $G$-equivariant semi-flow $\phi_t$ on $T_1\times T_2$ which restricts to the identity on  $\calc_\calq$,
and such that for all $x\in T_1\times T_2$, there exists $t\in [0,+\infty[$ such that $\phi_t(x)\in\calc_\calq$.

Moreover, for any rectangle $R$ such that for both $i\in\{1,2\}$, $p_i(R\cap\calc_\calq)=p_i(R)$,
$R$ is invariant under $\phi_t$.
\end{lem}

\begin{proof}[Proof of the proposition from the lemma]
Minimality and the connectedness of $\Hat\calc$ ensures that  $p_i(\Hat\calc)=T_i$.
Thus the hypotheses of the lemma are satisfied, so $\phi_t$ provides the required equivariant retraction by deformation
of $T_1\times T_2$ onto $\Hat\calc$.

The only thing remaining to check is that the trace of $\Hat\calc$ on a rectangle $R$ is either empty or contractible.
But if $R\cap \Hat\calc$ is nonempty, let $R_0=p_1(R\cap \Hat\calc)\times p_2(R\cap \Hat\calc)$ be the smallest rectangle containing $R\cap \Hat\calc$.
The lemma claims that $R_0$ is invariant under the semi-flow so $R\cap\Hat\calc$ is a retract by deformation of $R_0$.
\end{proof}

\begin{proof}[Proof of the lemma.]
For all $t\in [0,\infty]$ and $x\in T_1\times T_2$, we want to define $\phi_t(x)$ (see figure \ref{fig_vectorfield}).
If $x\in \calc_\calq$, we let $\phi_t(x)=x$ for all $t\geq 0$. 

Otherwise, using the fact that $p_i(\calc_\calq)=T_i$, 
choose  $a,b \in \calc$ such that $p_1(a)=p_1(x)$ and $p_2(b)=p_2(x)$ and let $R=[x_1,b_1]\times[x_2,a_2]$ be the smallest rectangle
containing $a$ and $b$ (we use the notation $a_i=p_i(a)$ and $b_i=p_i(b)$).
By proposition \ref{prop_connexity_coherent}, $\calc_\calq\cap R$ is connected. Moreover, $R\setminus\calc_\calq$ has at most
two connected components: the component $R_x$ containing $x$, and the component containing
the opposite corner $(b_1,a_2)$ if $(b_1,a_2)\notin \calc_\calq$ (indeed, the two components are the union of the traces of quadrants
containing the corners $x$ and $(b_1,a_2)$ respectively).
We say that a rectangle $R$ is a chart for $x$ if $x$ is a corner of $R$, and if
there are two points $a,b\in R\cap\calc_\calq$ having the same horizontal and vertical projection as $x$ respectively.

Given a chart $R$ for $x$, we identify it with  $[0,l_1]\times[0,l_2]\subset \bbR^2$
by sending $x$, $a$ and $b$ to $(0,0)$, $(l_1,0)$, $(0,l_2)$ where $l_1=d(x_1,b_1)$ and $l_2=d(x_2,a_2)$. 
Flow lines will be parallel to the vector $\vec v=(1,1)$.
Using this identification, given $y\in R$ and $t\in [0,\min(l_1,l_2)]$ not too large, it makes sense to write $y+t\vec v$.

Since $R\cap \calc_\calq$ is connected, there exists $t\leq\min(l_1,l_2)$ such that $x+t\vec v\in\calc_\calq$ since otherwise,
this segment would separate $a$ from $b$, contradicting the connectedness of $\calc_\calq\cap R$.
Similarly, for all $y\in R_x$, there exists $t\in [0,\min(l_1,l_2)]$ such that $y+t\vec v\in\calc_\calq$.
Thus, given the choice of $R$, we can define $\tau_R(y)$ to be the smallest positive $s$ such that $y+s\vec v\in\calc_\calq$.
We now define the flow on $R$ by $\phi^R_t(y)=y+\min(t,\tau_R(y))\vec v$.

The definition of the semi-flow does not change if we change a chart $R$ to a smaller one $R'$.
Indeed, the point $y+\tau_{R}(y)\vec v$ also lies $R'$ since otherwise, there would be
no $s$ such that $y+s\vec v\in \calc_\calq\cap R'$. This means that the definition of $\tau_R$ and $\tau_{R'}$
agree. Therefore, the definitions of $\phi_t^R$ and $\phi_t^{R'}$ agree.
Since there is a smallest chart for defining the flow at a given point $x$,
the definition of the semi-flow does not depend on any choice.

Note that by definition, $\phi_t(x)$ stays in any chart for $x$. Now if $R$ is a rectangle such that 
 for both $i\in\{1,2\}$, $p_i(R\cap\calc_\calq)=p_i(R)$, then $R$ contains a chart for each point of $R$.
In particular, $R$ is invariant under $\phi_t$. 
\\

Finally, we prove the continuity of $\phi:\bbR_+\times T_1\times T_2\ra T_1\times T_2$ by proving 
that the semi-flow is Lipschitz with respect to each variable separately. On $T_1\times T_2$, we consider the distance
$d(x,y)=max_{i\in\{1,2\}} d_{T_i}(x_i,y_i)$.
Clearly, $t\mapsto\phi_t(x)$ is $1$-Lipschitz.
Now consider $x_1\mapsto \phi_t(x_1,x_2)$ where $t$ and $x_2$ are fixed. 
We choose $b_1\in T_1$ such that $b=(b_1,x_2)\in \calc_\calq$.
Let $x_1,x'_1\in T_1$, and denote $x=(x_1,x_2)$ and $x'=(x'_1,x_2)$.
Let $l=d(x_1,x'_1)$.

\begin{figure}[htbp]
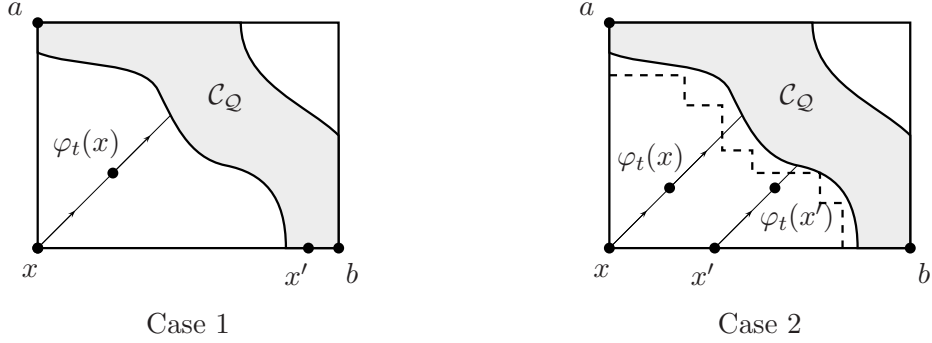

  \centering
  \input flot_lipschitz.pst
  \caption{The flow is Lipschitz}
  \label{fig_flot_lipschitz}
\end{figure}

\paragraph{Case 1:} $x \notin\calc_\calq$ and $x' \in\calc_\calq$.
In this case, choose $R=[x_1,x'_1]\times [x_2,a_2]$ as a chart, where $(x_1,a_2)\in \calc_\calq$.
Then, $d(\phi_t(x),\phi_t(x'))=d(\phi_t(x),x')\leq l$.
Indeed let $\lambda=\min(t,\tau(x))$,
and note that $\lambda\leq l$ because $\phi_t(x)\in R$.
Then $d(\phi_t(x),x')= ||\lambda\vec v-x'||=\max(|\lambda|,|\lambda-l|)\leq l$.

\paragraph{Case 2:} $x$ and $x'$ don't lie in $\calc_\calq$,
and $x'_1\in [x_1,b_1]$ (remember that $b_1$ is such that $(b_1,x_2)\in\calc_\calq$). 
In this case, we choose $R=[x_1,b_1]\times [x_2,a_2]$
as a chart where $a=(x_1,a_2)\in\calc_\calq$. The arc $[x_1,x'_1]\times \{x_2\}$ cannot
intersect $\calc_\calq$ because fibers are convex and $(b_1,x_2)\in \calc_\calq$.
Thus $x$ and $x'$ are in the same component $R_x$ of $R\setminus \calc_\calq$.
Let $\lambda=\min(t,\tau(x))$ and $\lambda'=\min(t,\tau(x'))$.
Then  $d(\phi_t(x),\phi_t(x'))=||(\lambda,\lambda)-(\lambda'+l,\lambda')||\leq l+|\lambda-\lambda'|$.
There remains to prove that $z_1 \mapsto\tau(z_1,x_2)$ is a Lipschitz function in restriction to $R_x$.
Write $R_x$ as an increasing union of sets $R_k$ where each 
$R_k$ is a finite union of traces of light quadrants on $R$.
It is easy to check that the piecewise linear map $\tau_k$ obtained by
replacing $R_x$ by $R_k$ is $\sqrt{2}$-Lipschitz. 
The function $\tau$ being the supremum of the fuctions $\tau_k$,
it is therefore $\sqrt{2}$-Lipstchitz
\\

We now show how to deduce the other cases from case 1 and 2:
we just have to consider the case where $x,x'\notin \calc_\calq$.
Consider the fiber $\calf_1=\{b_1\in T_1\,|\, (b_1,x_2)\in\calc_\calq\}$. 
If $\calf_1$ meets $[x_1,x'_1]$ in a point $b'_1$,
we apply case 1 to the pair $x,b'$ and to the pair $b',x'$ and apply triangle inequality.
If $\calf_1$ does not meet $[x_1,x'_1]$, let $x''_1$ be the center of the tripod $(x_1,x'_1,b_1)$.
Now apply case 2 to the pair $x,x''$ and to the pair $x'',x'$.
\end{proof}

\subsection{The augmented core is CAT(0).}

\newcommand{\length}{{\mathrm{length}}}

Let $d$ be the usual CAT(0) metric $d$ on $T_1\times T_2$ defined by
$$d((a_1,a_2),(b_1,b_2))=\sqrt{d_{T_1}(a_1,b_1)^2+d_{T_2}(a_2,b_2)^2}.$$
Let $d_{\Hat\calc}$ be the path metric  induced by $d$ on $\Hat\calc$
\[ d_{\Hat\calc}(x,y)=\inf \{\length(c)|\ c\text{ path joining $x$ to $y$ in $\Hat\calc$}\} \]
where the length is measured with $d$.

\begin{prop}\label{prop_cat0}
 When endowed with the metric $d_{\Hat\calc}$ above, the augmented core $\Hat\calc$ is CAT(0).
\end{prop}

\begin{proof}
First, the result is clear when $T_1$ and $T_2$ are simplicial trees since
$\Hat\calc$ is a simply connected square complex, and the link at each vertex has
no  has no non-trivial loop of length less than $2\pi$
because this is already true in $T_1\times T_2$. 

 For $\bbR$-trees, first note that the metric  $d_{\Hat\calc}$ is well defined and finite
since any two points are joined by a Lipschitz path (\ref{prop_contractible}) (note however that it might not be complete as
$T_i$ is itself usually not complete).
Then, if $K_1,K_2$ are convex subsets of $T_1$ and $T_2$ respectively,
for any Lipschitz path $c:[0,1]\ra \Hat\calc$ joining two points of $K_1\times K_2\cap \Hat\calc$
there is a shorter path in  $K_1\times K_2\cap \Hat\calc$ joining them:
for each component $]a,b[$ of $[0,1]\setminus (p_1\circ c)\m(K_1)$, 
one has that $p_1(c(a))=p_1(c(b))\in K_1$, and one can replace $c$ on $[a,b]$ by the 
geodesic in the vertical fiber of $p_1(c(a))$ and get a shorter path in $\Hat\calc$ with the same Lipschitz constant.
Doing this for every such component, we get a shorter path $c'$ such that $p_1\circ c'([0,1])\subset K_1$.
Doing the symmetrical operation on $c'$, we get a shorter path in $K_1\times K_2\cap \Hat\calc$.
Given $a=(a_1,a_2),b=(b_1,b_2)\in \Hat\calc$, one can take $K_i$ to be the compact interval $K_i=[a_i,b_i]$,
and the compactness of $K_1\times K_2\cap \Hat\calc$ implies that the infimum 
in $d(a,b)$ $\Hat\calc$ is achieved. 
In particular $\Hat\calc$ is geodesic, and sets of the form
$K_1\times K_2\cap\Hat \calc$ are convex.

Moreover, there is a unique geodesic between two given points. Otherwise, one can find two points $a,b$ with
two geodesics $c_1$, $c_2$ joining them with $c_1\cap c_2=\{a,b\}$. The fact that the trace of a rectangle is convex in $\Hat\calc$ 
implies that $c_1$ and $c_2$ are graphs of monotonous functions in the smallest rectangle containing $\{a,b\}$.
We now see $c_1$, $c_2$ as maps of monotonous functions $I_1\ra I_2$, and assume for instance that $c_1\leq c_2$. 
Since the set lying between $c_1$ and $c_2$ is contained in $\Hat\calc$,
if $c_1$ is not concave, then one could shorten it. Similarly, $c_2$ is necessarily convex.
This implies $c_1=c_2$.

Now to prove that $\Hat\calc$ is CAT(0), since any geodesic triangle is contained in the product of two tripods,
we need only to prove that $K_1\times K_2\cap \Hat\calc$ is CAT(0) where $K_i$ is a (compact) tripod in $T_i$. 
Let $X=\Hat\calc\cap K_1\times K_2$. Since $X$ can be obtained by removing from $K_1\times K_2$
a countable set of quadrants, write $X$ as a decreasing intersection of a sequence of sets $X_k$
obtained by removing finitely many quadrants from $K_1\times X_2$.
The argument in the simplicial setting implies that $X_k$ is CAT(0).
Consider a sequence of linearly parametrized geodesics $c_k:[0,1]\ra X_k$ joining two points $a,b\in X$.
Up to extracting a subsequence, $c_k$ converges to a curve $c$ joining $a$ to $b$ in $X$.
General nonsense shows that $\length(c)\leq \lim \length(c_k)$ and since $X\subset X_k$, $\length(c)\geq \length(c_k)$,
so the path metrics on $X_k$ converge to the path metric on $X$.
It follows that $X$ is CAT(0).
\end{proof}

\section{Characterization of the core}\label{sec_characterization}

\begin{prop}[Characterization of the core]\label{prop_carac_core}
Let $T_1,T_2$ be two actions of a group $G$ on $\bbR$-trees such that $\calc\neq\es$.
  Let $F\subset T_1\times T_2$ be a non-empty closed connected $G$-invariant subset
with convex fibers.
Then $F$ contains $\calc(T_1\times T_2)$. 

Moreover, $\calc$ is the intersection of all such sets $F$.
\end{prop}

If $T_1$ and $T_2$ are not the refinements of a common splitting, then $\calc$ is itself a closed connected subset with convex fibers.
We thus get:

\begin{SauveCompteurs}{cor_carac_core}
\begin{cor}\label{cor_carac_core}
Let $T_1$, $T_2$ be two actions of $G$ on $\bbR$-trees whose minimal subtrees are dense. 
Assume that $\calc$ is non-empty and that $T_1$ and $T_2$ are not the refinement of a common simplicial non-trivial action.

Then $\calc$ is the smallest non-empty closed invariant connected subset of $T_1\times T_2$
having convex fibers.
\end{cor}
\end{SauveCompteurs}

We will often use this characterization of the core under the following form:
\begin{cor}\label{cor_maps}
Let $T_1,T_2$ be two actions of a group $G$ on $\bbR$-trees.
  Let $X$ be a nonempty connected space with an action of $G$ such that there are two equivariant maps
$f_1$, $f_2$ from $X$ to $T_1$ and $T_2$ such that the preimage of each point of $T_i$ is connected.
Let $F=(f_1,f_2):X\ra T_1\times T_2$.

Then $\ol{F(X)}$ contains $\calc$.
\end{cor}

\begin{proof}[Proof of corollary \ref{cor_maps}]
 Since $\ol{F(X)}$ is an nonempty invariant closed connected subset of $G$, we just have to prove that
 it has connected fibers. By Corollary \ref{cor_adherence} below, we need only to check that $F(X)$ has connected fibers.
So let $x\in T_1$, and consider a fiber $F(X)\cap p_1\m(x_1)$.
But this fiber can be written as $F(f_1\m(x_1))$, which is connected since $f_1\m(x_1)$ is connected.
\end{proof}

The main result to prove the proposition is the following lemma.

\begin{lem}\label{lem_fibre_quadrant}
Let $T_1$, $T_2$ be two $\bbR$-trees and let $F$ be a nonempty connected subset of $T_1\times T_2$
with convex fibers.
Then the complement of $\ol F$ is a union of quadrants.
\end{lem}

This corollary follows immediately:
\begin{cor}\label{cor_adherence}
  If $F$ is nonempty, connected and has convex fibers, then so is $\ol F$.
\end{cor}

\begin{rem*}
  Note that this is of course false if one removes any connectedness assumption:
just take $E=\{(x,x)|\, x\in\bbQ\} \cup \{(x,-x)|\, x\in \bbR\setminus\bbQ\}$ in $\bbR^2$.
\end{rem*}

\begin{proof}[Proof of the proposition from the lemma]
  Let $\calq$ be the family of quadrants which don't intersect $F$ so that, 
by lemma \ref{lem_fibre_quadrant}, $F=T_1\times T_2 \setminus \bigcup_{Q\in\calq}Q$.
To prove that $F$ contains $\calc$, we only need to prove that any quadrant $Q\in\calq$
is light. By choosing the base point in $F$, the orbit of the base point does not meet any quadrant $Q\in\calq$
so the proposition follows.
\end{proof}

We will use the following terminology: if $Q=\delta_1\times\delta_2$ is a quadrant based at $x=(x_1,x_2)$,
we call $\partial_1 Q=\{x_1\}\times \delta_2$ (resp. $\partial_2 Q=\delta_1\times\{x_2\}$)
the \emph{vertical} (resp. the \emph{horizontal}) boundary of $Q$. 
Note that $\partial Q=\partial_1 Q\disjoint \partial_2 Q \disjoint \{x\}$.
Say that two quadrants $\delta_1\times\delta_2$ and $\delta'_1\times\delta'_2$ based at the same point are \emph{opposite}
if $\delta_1\neq\delta'_1$ and $\delta_2\neq\delta'_2$.

Let's start with the following fact:

\begin{fact}\label{fact_quadrants}
Let $T_1\times T_2$ be two $\bbR$-trees and let $F$ be a nonempty connected subset of $T_1\times T_2$
with convex fibers.
Fix a point $x\notin F$.
  \begin{enumerate*}
  \item \label{exist_dble} If $F$ meets two opposite quadrants based at $x$, then there exists a quadrant $P$,
based at $x$, such that $F$ intersects the horizontal and the vertical boundary of $P$.
\item\label{un_seul_dble}  Let $P$ be a quadrant based at $x$ such that $F$ intersects the horizontal and the vertical boundary of $P$.
Then the closure of any quadrant opposite to $P$ doesn't intersect $F$.
In other words, any quadrant whose closure intersects $F$ has a common (vertical or horizontal) boundary with $P$.
Moreover, $F$ meets $P$.
  \end{enumerate*}
\end{fact}

\begin{proof}[Proof of the fact]
\ref{exist_dble}.
Let $Q=\delta_1\times\delta_2$, $Q'=\delta'_1\times\delta'_2$ be two opposite quadrants based at $x$ which intersect $F$.
Since $F$ is connected and does not contain $x$, $F$ meets $\partial_1 Q$ or $\partial_2 Q$.
Assume for instance that  $F$ meets $\partial_1 Q=\{x_1\}\times\delta_2$.
For the same reasons, $F$ meets $\partial_1 Q'$ or $\partial_2 Q'$.
Since fibers of $F$ are convex, and since $x\notin F$, $F$ cannot intersect $\partial_1 Q'$.
Therefore $F$ meets both boundaries of the quadrant $P=\delta'_1\times\delta_2$. 

\ref{un_seul_dble}.
Now let $P=\rho_1\times\rho_2$ be a quadrant based at $(x_1,x_2)$ such that $F$ intersects $\partial_1 P$ and $\partial_2 P$.
Convexity of fibers implies that $F$ cannot intersect the set $A=\rho_1^*\times\{x_2\}\cup \{x_1\}\times \rho_2^*$.
Now let $Q$ be a quadrant opposite to $P$. Since $\partial Q\subset A\cup\{x\}$, 
the connectedness of $F$ prevents $F$ from intersecting $\ol Q$.
To prove that $F$ meets $P$, just note that $A\cup P\cup \{x\}$ separates $\partial_1 P$ from $\partial_2 P$.
\end{proof}

\begin{proof}[Proof of Lemma \ref{lem_fibre_quadrant}]
  Let $x=(x_1,x_2)\notin \ol F$. We have to find a quadrant $Q$ containing $x$ and disjoint from $F$.
Let $V=V_1\times V_2$ be an open neighbourhood of $x$ which does not intersect $F$.
If $F$ does not intersect any quadrant based at $x$, then for instance, one can assume that $F$ is contained
in $\{x_1\}\times \delta_2$ for some direction $\delta_2$ at $x_2$, and the result is clear.

Let $a\in F$ lying in a quadrant based at $x$. Let $y\in V\cap \, ]x_1,a_1[\times]x_2,a_2[$,
and let $Q_y$ be the quadrant based at $y$ containing $x$ (see figure \ref{fig_cas_pres} and \ref{fig_cas_loin}). 
If $Q_y\cap F=\es$, then we are done.
Otherwise, $F$ meets two opposite quadrants based at $y$ (namely $Q_y$ and the quadrant containing $a$).
The fact above says that there is a quadrant $P_y=\rho_1\times\rho_2$ based at $y$ such that $F$ meets both the horizontal
and the vertical boundary of $P_y$. Since $P_y$ has a common boundary with every quadrant at $y$ meeting $F$,
we can assume for instance that $P_y$ and $Q_y$ share their vertical boundary.
Consider a point $v=(y_1,v_2)\in F\cap \partial_1 P_y$, and a point $u=(u_1,y_2)\in F\cap\partial_2 P_y$.
We distinguish two cases:
\begin{enumerate*}
\item\label{cas_pres} either $x_2\notin [v_2,y_2]$
\item \label{cas_loin} or $x_2\in [v_2,y_2]$
\end{enumerate*}

First assume that case \ref{cas_pres} occurs (figure \ref{fig_cas_pres}). 
Since $x\in Q_y$, one has that $x_2$ and $v_2$ are in a common direction based at $y_2$ (namely, the vertical direction of $Q_y$).
Therefore, $[y_2,v_2]$ and $[y_2,x_2]$ have a common nondegenerate initial segment $[y_2,z_2]$. 
One has $z_2\neq x_2$ since otherwise case \ref{cas_loin} would occur, and $z_2\neq v_2$ because $F$ does not intersect $V$.
This means that the three directions $\delta(y_2)$, $\delta(x_2)$ and $\delta(v_2)$ 
based at $z_2$ and containing respectively $y_2$, $x_2$, and $v_2$ are distinct.
Since $x_2,z_2,y_2,a_2$ are aligned in this order, $a_2\in\delta(y_2)$.

\begin{figure}[htbp]
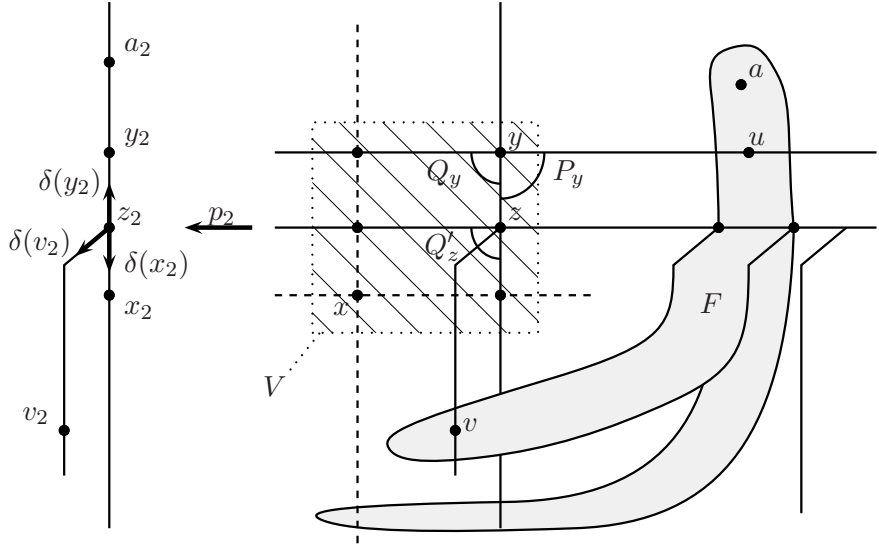

  \centering
  \input cas_pres.pst
  \caption{Proof of lemma \ref{lem_fibre_quadrant}, case \ref{cas_pres}}
  \label{fig_cas_pres}
\end{figure}

Now consider the point $z=(y_1,z_2)\in V$. In particular $z\notin F$.
Let $Q'_{z}=\delta_1\times\delta(x_2)$ be the quadrant at $z$ containing $x$.
If $Q'_{z}$ does not intersect $F$, we are done. Otherwise, there are two opposite quadrants at $z$ meeting $F$
(namely $Q'_{z}$ and the quadrant $\eta_1\times\delta(y_2)$ containing $a$). The quadrant $P'_{z}$ given by point \ref{exist_dble} 
of the fact is one of the two quadrants $\delta_1\times\delta(y_2)$ or $\eta_1\times\delta(x_2)$.
Since $F$ intersects the vertical boundary of $P'_{z}$, $F$ contains a point in $\{y_1\}\times\delta(y_2)$
or in $\{y_1\}\times\delta(x_2)$. Since $v \in \{y_1\}\times \delta(v_2)$ also lies in $F$, convexity of fibers
contradicts the fact that $z\notin F$.
This concludes case \ref{cas_pres}.\\

\begin{figure}[htbp]
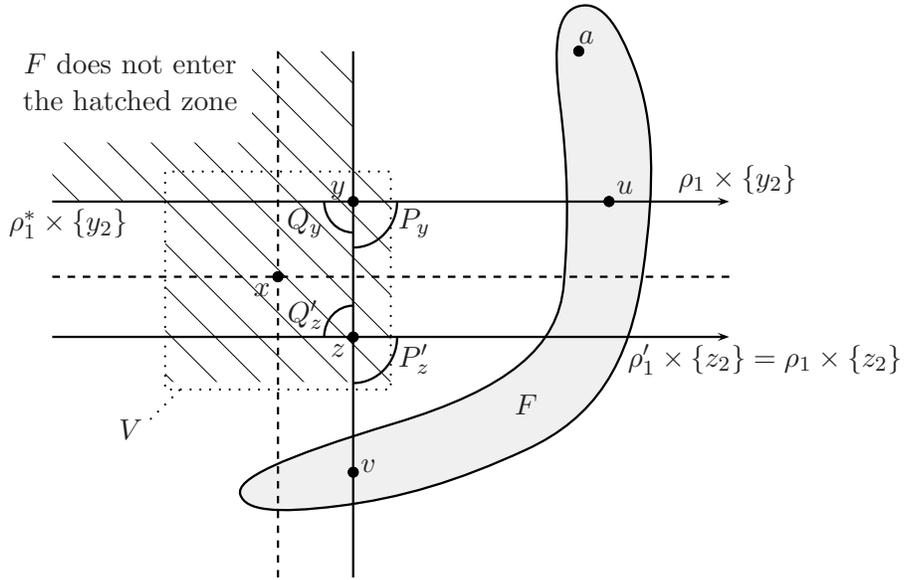

  \centering
  \input cas_loin.pst
  \caption{Proof of lemma \ref{lem_fibre_quadrant}, case \ref{cas_loin}}
  \label{fig_cas_loin}
\end{figure}

Now we assume that $x_2\in [v_2,y_2]$ (see figure \ref{cas_loin}).
Let $z=(y_1,z_2)$ where $z_2\in V_2\cap \, ]x_2,v_2[$.
We are going to prove that the quadrant $Q'_z$, based at $z$ and containing $x$
does not meet $F$ by finding a quadrant $P'_z$ opposite to $Q'_z$, and meeting $F$ at its two boundaries.
The fact will then conclude.

Let $\rho'_2$ be the direction at $z_2$ containing $v_2$, and let $\rho'_1=\rho_1$.
Since we are in case \ref{cas_loin}, $P'_z=\rho'_1\times\rho'_2$ is opposite to $Q'_z$.
Moreover, the choice of $\rho'_2$ implies that $F$ meets the vertical boundary of $P'$.
Now notice that $\rho_1^*\times\{y_2\}\cup V \cup \rho_1\times\{z_2\}$ separates $u$ from $v$.
Moreover, $\rho_1^*\times\{y_2\}$ does not intersect $F$ since it is contained in the closure
of the union of quadrants opposite to $P_y$. Nor does $V$ intersect $F$. Thus, since $F$ contains $u$ and $v$,
$F$ must intersect $\rho_1\times\{z_2\}$, which is the horizontal boundary of $P'_z$.
Thus $F$ meets both boundaries of $P'_z$, which implies that $F$ does not meet $Q'_z$,
so $Q'_z$ is the desired quadrant.
\end{proof}

\section{Compatibility of tree actions}\label{sec_compatibility}

We say that $\calc$ is \emph{1-dimensional} if it does not contain any rectangle $I_1\times I_2$ 
where $I_1$ and $I_2$ are nondegenerate arcs. The proof will show that when $\calc$ is 1-dimensional,
$\Hat\calc$ has a natural structure of $\bbR$-tree, and is a common refinement of $T_1$ and $T_2$.

\begin{SauveCompteurs}{cpt_compatible}
\begin{thm}[compare \cite{ScSw_splittings}]\label{thm_compatible}
   Let $T_1,T_2$ be two minimal actions of $G$ on $\bbR$-trees, such that $\calc(T_1\times T_2)\neq \es$.

Then $T_1$ and $T_2$ have a common refinement if and only if $\calc$ is $1$-dimensional.
\end{thm}
\end{SauveCompteurs}

\begin{rem*}
  The result does not hold if we don't assume that $\calc\neq\es$:
an abelian action having no invariant line, and the corresponding action by translation on $\bbR$
have no common refinement. 
\end{rem*}

\begin{proof}
First assume that $T_1$ and $T_2$ have a common refinement $T_0$,
and denote by $f_i:T_0\ra T_i$ an equivariant map preserving alignement.
Consider the map $F=(f_1,f_2):T_0\ra T_1\times T_2$.
By corollary \ref{cor_maps}, $\ol{F(T_0)}$ contains $\calc$.

Now, $F(T_0)$ is $1$-dimensional: this is clear in the simplicial context; 
in general, this can be proved as follows. Take
$a=(a_1,a_2)$, $b=(a_1,c_2)$, $c=(c_1,c_2)$, $d=(c_1,a_2)$ 
be the corners of a rectangle contained in $F(T_0)$, and let $a_0$, $b_0$, $c_0$, $d_0$
be some preimages in $T_0$. Let $K_0\subset T_0$ be the convex hull of $\{a_0,b_0,c_0,d_0\}$.
Then one has either $[a_0,b_0]\cup [c_0,d_0]=K_0$ or $[b_0,c_0]\cup [d_0,a_0]=K_0$.
Let 's assume for instance that  $[a_0,b_0]\cup [c_0,d_0]=K_0$. Now since $f_1(a_0)=f_1(b_0)$,
$f_1$ is constant on $[a_0,b_0]$. For similar reasons, $f_1$ is constant on $[c_0,d_0]$.
Thus $f_1$ is constant on $K_0$, a contradiction.
One deduces that $\ol{F(T_0)}$ is also $1$-dimensional: consider a non-degenerate rectangle $R\subset \ol{F(T_0)}$
and let $x\in R$ which does not lie on its boundary. Then the four quadrants based at $x$ and 
containing the four corners of $R$ intersect $F(T_0)$. By fact \ref{fact_quadrants},
 $x\in F(t_0)$.\\

Conversely, assume that $\calc$ and thus $\Hat\calc$ is one dimensional. 
If $T_1$ and $T_2$ are simplicial trees, then $\Hat\calc$ is a contractible one-complex, \ie a simplicial tree.
In general, we will prove that $\Hat\calc$ is an $\bbR$-tree when endowed with the metric
$d((x_1,x_2),(y_1,y_2))=d_{T_1}(x_1,y_1)+d_{T_2}(x_2,y_2)$.
One can deduce that $\Hat\calc$ itself is a common refinement for $T_1$ and $T_2$:
the map $p_i:\Hat\calc\ra T_i$ preserves alignement
because $\Hat\calc$ has convex fibers.

Let's prove that $\Hat\calc$ is an $\bbR$-tree. 
We have already proved that $\Hat\calc$ is CAT(0) and therefore geodesic.
Assume that $\Hat\calc$ contain an embedded circle $c=(c_1,c_2):S^1\ra\Hat\calc$.
Then $c_1$ is necessarily non constant, so there exists a non-degenerate interval $[a,b]\subset S^1$
with $c_1(a)=c_1(b)\notin c_1(]a,b[)$. 
Note that the segment $\{c_1(a)\}\times [c_2(a),c_2(b)]$ is contained in $\Hat\calc$ by convexity of fibers.
Now $c_1(]a,b[)$ is contained in a direction based at $c_1(a)$,
so for all $\eps>0$, there exists $a'\in(]a,b[)$ close to $a$ and $b'\in]a,b[$ close to $b$ 
such that $c_1(a')=c_1(b')$, $d(c_2(a),c_2(a'))<\eps$ and $d(c_2(b),c_2(b'))<\eps$.
If $\eps$ was chosen very small compared to $d(c_2(a),c_2(b))$, then the segments $[c_2(a),c_2(b)]$ and $[c_2(a'),c_2(b')]$
in $T_2$ must intersect in a nondegenerate interval, say $I_2$.
Now the segments $\{c_1(a)\}\times I_2$ and $\{c_1(a')\}\times I_2$ are contained in $\Hat\calc$ by convexity of fibers,
which implies that the rectangle $[c_1(a),c_1(a')]\times I_2$ is contained in $\Hat \calc$, a contradiction.

\end{proof}

\section{Topological interpretation of the intersection number}\label{sec_interpretation}
We now give a topological interpretation of the intersection number of two splittings.

We need a few definitions. Given a connected cell complex $X$ and a (maybe disconnected) subcomplex $Y\subset X$,
say that $Y$ is \emph{$2$-sided} if $Y$ has a neighbourhood in $X$ homeomorphic to  $Y\times[-\eps,\eps]$.
Let $\Tilde X$ be its universal covering, and $\Tilde Y$ be the preimage of $Y$ in $\Tilde X$.
The \emph{tree dual to $Y\subset X$} is the graph $T_Y$ whose vertices are connected components of $\Tilde X\setminus \Tilde Y$
and whose edges are connected components of $\Tilde Y$, and an edge $e$ is adjacent to a vertex $v$ if $e\subset \ol{v}$. The simple connectedness 
of $\Tilde X$ implies that $T_Y$ is a tree, and it is clearly endowed with an action of $\pi_1(\Sigma)$.
Given $Y_1,Y_2$ two $2$-sided subcomplexes of $X$, say that $Y_1$ and $Y_2$ intersect transversely if 
$Y_1\cap Y_2$ has a neighbourhood $N$  homeomorphic to $Y_1\cap Y_2\times [-\eps,\eps]^2$
where $Y_i\cap N$ is mapped to $Y_1\cap Y_2\times p_i\m(0)$ where $p_i:[-\eps,\eps]^2\ra [-\eps,\eps]$ is the canonical projection.

In the following proposition, $\#\pi_0(Y_1\cap Y_2)$ denotes the number of connected components of $Y_1\cap Y_2$.

\begin{SauveCompteurs}{cpt_interpretation}
\begin{thm}\label{thm_interpretation}
  Assume that $Y_1,Y_2\subset X$ are two 2-sided subcomplexes, which intersect transversely, and let $T_1$, $T_2$ be the two dual trees, endowed with
the action of $\pi_1(X)$. Then $i(T_1,T_2)\leq \#\pi_0(Y_1\cap Y_2)$.

Moreover, given two non-trivial actions of a group on simplicial trees $T_1,T_2$,
there exists a complex $X$ and $Y_1,Y_2\subset X$ two 2-sided subcomplexes intersecting transversely such that 
 $i(T_1,T_2)= \#\pi_0(Y_1\cap Y_2)$.
\end{thm}
\end{SauveCompteurs}

\begin{proof}
Assume that $T_Y$ is dual to $Y\subset X$.
Then there is a continuous equivariant map $f:\Tilde X\ra T_Y$ defined by sending each connected component of $\Tilde Y\times ]-\eps,\eps[$
to an edge, and by sending each connected component of $\Tilde X\setminus \Tilde Y\times ]-\eps,\eps[$ to a vertex.
The main observation here is that the preimage of each point of $T_Y$ under $f$ is connected.

In our setting, denote by $f_i:\Tilde X\ra T_i$ the maps defined above, and let $F=(f_1,f_2)$.
By corollary \ref{cor_maps}, $F(\Tilde X)$ contains $\calc$ since $F(\Tilde X)$ is closed as a subcomplex of $T_1\times T_2$.
Since $Y_1$ and $Y_2$ intersect transversely, $Y_1\cap Y_2$ has a neighbourhood homeomorphic to $Y_1\cap Y_2\times [-\eps,\eps]^2$,
so $F$ maps each connected component of $\Tilde Y_1\cap \Tilde Y_2\times [-\eps,\eps]$ to a $2$-cell of $T_1\times T_2$.
Moreover, the preimage of the center of a $2$-cell of $T_1\times T_2$ is a connected component of 
$\Tilde Y_1\cap \Tilde Y_2$.
Therefore, the number of orbits of $2$-cells in $F(\Tilde X)$ is bounded by the number of orbits of connected components of $\Tilde Y_1\cap \Tilde Y_2$
so 
$$\begin{array}{rcccl}
i(T_1,T_2)&=&\#\{\text{$2$-cells of } \calc/G\}
&\leq & \#\{\text{$2$-cells of } F(\Tilde X)/G\}\\
&\leq& \#\  \pi_0(\Tilde Y_1\cap \Tilde Y_2)/G &=& \#\pi_0(Y_1\cap Y_2).
\end{array}$$

We now prove that equality can achieved by constructing $X,Y_1,Y_2$ from $\Hat\calc$. 
Let $\Tilde X=\Hat\calc$, let $E_i$ be the set of midpoints of edges in $T_i$,
and let $\Tilde Y_i=\Hat\calc\cap p_i\m(E_i)$. In the category of complexes of groups, one could take $Y_i=\Tilde Y_i/G$
and $X=\Tilde X/G$. However, we need to modify this construction to get free actions.
Take $A$ a simply connected complex on which $G$ acts freely (for instance, a Cayley $2$-complex).
Let $\Tilde X'=\Hat\calc\times A$ endowed with the diagonal action of $G$,
and let $\Tilde Y'_i=\Tilde Y_i\times A$. 
Connected components of $\Tilde Y'_1\cap \Tilde Y'_2$ are of the form $x\times A$ where $x$ is either the center of a $2$-cell of $\calc$,
or the midpoint of the main diagonal of a maximal twice light rectangle in $\Hat\calc$. 
Thus, in the presence of twice light rectangles, we need to change $E_i$ so that $\Tilde Y'_1\cap \Tilde Y'_2$
contains no point in the main diagonal of a twice light rectangle. To this means, one can keep $E_1$ unchanged, and take for $E_2$ 
an equivariant set of points meeting each edge of $T_2$ exactly once, but not containing the midpoint of any edge.
With this modification, there is a one-to-one correspondance between connected components of $\Tilde Y'_1\cap\Tilde Y'_2$
and the $2$-cells of $\calc$. In particular, $\# \pi_0(\Tilde Y'_1\cap \Tilde Y'_2)/G= i(T_1,T_2)$.

Let $X'=\Tilde X'/G$ and $Y'_i=\Tilde Y'_i/G$.
Since $\Tilde X'$ is simply connected and since the action of $G$ on $\Tilde X'$ is free,
$\Tilde X'$ is the universal cover of $X'$ and $G\simeq\pi_1(\Sigma)$.
Thus, $T_i$ is dual to $Y'_i\subset X'$, and $i(T_1,T_2)=\#\pi_0(Y'_1\cap Y'_2)$.
\end{proof}

\section{Core of geometric actions}\label{sec_geometric}

The goal of this section is to produce a finite \emph{fundamental domain} for the core of geometric actions in the 
following weak sense:

\begin{SauveCompteurs}{cpt_fund_domain}
\begin{thm}\label{thm_fund_domain}
  Let $T_1$, $T_2$ be geometric actions of a finitely generated group $G$ on $\bbR$-trees.

Then there is a set $D\subset T_1\times T_2$, which is a finite union of compact rectangles,
and such that $\calc(T_1\times T_2) \subset\ol{G.D}$.
\end{thm}
\end{SauveCompteurs}

Remember that a minimal action of a finitely generated group on a simplicial tree is geometric
if and only if its edge stabilizers are finitely generated. Therefore, we get:

\begin{cor}[\cite{Sco_symmetry}]\label{cor_finite_intersection_simplicial}
  Let $T_1,T_2$ be two splittings of a finitely generated group $G$ over finitely generated groups.

Then $\calc(T_1\times T_2)/G$ is compact. In particular, $i(T_1,T_2)$ is finite.
\end{cor}

In general, we will need a stronger assumption to deduce the finiteness of the intersection number:
\begin{prop}\label{prop_finite_intersection}
 Let $T_1$, $T_2$ be geometric actions of a \emph{finitely presented} group $G$ on $\bbR$-trees.

Then $i(T_1,T_2)$ is finite.
\end{prop}

The philosophy here is the following: we construct $2$-complex $X$, with a cocompact action of $G$,
 with two measured foliations $\calf_1,\calf_2$ so that $T_1$ and $T_2$
are the leaf spaces made Hausdorff of those foliations. Now let $f_i:X\ra T_i$ be the canonical projections and
let $F=(f_1,f_2):X\ra T_1\times T_2$.
First, $F(X)$ in contained in the orbit of a finitely many rectangles.
Now if points of $T_i$ exactly coincide with leaves of $\calf_i$, the fibers of $f_i$ are connected and Corollary \ref{cor_maps}
implies that $\ol{F(X)}$ contains $\calc$. 
In general, the connexity of fibers might fail, but we will get around this difficulty.

\subsection{An example of infinite intersection number}

To motivate this section, we first give an example of actions of a free group on simplicial trees $T_1,T_2$
such that $i(T_1,T_2)$ is infinite. This answers a question asked by Scott and Swarup in \cite{ScSw_splittings}.

\begin{lem}\label{lem_infinite_intersection}
Let $T_1$ be a \emph{free} minimal action of the free group $G=\langle a,b,c \rangle$ on a simplicial tree, for instance on its Cayley graph.
Let $H$ be a non-finitely generated subgroup of $\langle a,b \rangle$, and let $T_2$ be the Bass-Serre tree of the amalgam
$G=\langle a,b \rangle *_H (H*\langle c \rangle)$.

Then $i(T_1,T_2)=\infty$.
\end{lem}

\begin{proof}
First $\calc$ is nonempty and has no twice light rectangle; and by minimality $p_i(\calc)=T_i$;
now let $e_2$ be an edge of $T_2$ stabilized by $H$, and let $A=\calc\cap p_2\m(e_2)$.
Since fibers are convex, $A$ has the form $A_1\times e_2$ where $A_1$ is a nonempty subtree of $T_1$.
Note that the invariance of $\calc$ implies that $A_1$ is $H$-invariant, and
that two edges in $A_1$ are in the same $H$-orbit if and only if the corresponding rectangles
of $A$ are in the same $G$-orbit. Thus proving $i(T_1,T_2)=\infty$ consists in proving that $A_1$ has infinitely many $H$-orbits of edges.
But the action of $H$ on $A_1$ is free so $H$ occurs as the fundamental group of the graph $A_1/H$, which cannot be finite since
$H$ is not finitely generated.
\end{proof}

\subsection{Preliminaries about geometric actions}

All the material in this section is borrowed from \cite{LP} where more details can be found.
A \emph{measured foliation} $\calf$ on a $2$-complex $X$
consists of the choice, for each closed simplex $\sigma$ of $X$ of a (maybe constant) affine map
$f_\sigma:\sigma\ra \bbR$ defined up to post-composition by an isometry of $\bbR$, and such that 
is consistent under restriction to a face: if $\tau$ is a face of $\sigma$, then 
$f_\tau=\phi\circ(f_\sigma)_{|\tau}$ for some isometry of $\phi$ of $\bbR$.
Level sets of $f_\sigma$ give a foliation on each closed simplex. 
We say simply a \emph{foliated 2-complex} to mean a $2$-complex endowed with a measured foliation.
Leaves of the foliations on $X$ are defined as the equivalence classes of the equivalence relation generated 
by the relation \emph{$x,y$ belong to a same closed simplex $\sigma$ and $f_\sigma(x)=f_\sigma(y)$}.
The transverse measure $\mu(c)$ of a path $c:[0,1]\ra \sigma$ transverse (resp. parallel) to the foliation is the length of the interval $f_\sigma(c([0,1]))$. The transverse measure is invariant under the holonomy along the leaves.
The transverse measure thus defines a metric on each transverse edge.

We say that map $f$ from a simplex to an $\bbR$-tree $T$ is \emph{affine}
if $f=i\circ f_a$ where $a:\sigma\ra\bbR$ is an affine map,
and $i:I\ra T$ is an isometry defined on a convex subset $I\subset \bbR$ containing $a(\sigma)$.
Given a $2$-complex $X$ and a map $f:X\ra T$ which is affine in restriction to each simplex,
there is a natural measured foliation $\calf$ on $X$ defined by the restrictions of $f$ to the simplices of $X$.
We call $\calf$ the measured foliation induced by $f$.

If $c:[0,1]\ra X$ is a path wich is piecewise transverse or parallel to the foliation, we define $\mu(c)$ 
is the sum of the transverse measures of the pieces. 
The pseudo-metric 
$$\delta(x,y)=\inf\{\mu(c)\text{ for $c$ joining $x$ to $y$}\}$$
is zero on each leaf of $X$. 
By definition, the \emph{leaf space made Hausdorff} $X/\calf$ of $X$ is the metric space obtained from $X$
by making $\delta$ Hausdorff, \ie by identifying points at pseudo-distance $0$. 

\begin{thm}[\cite{LP}] \label{thm_LP_leaves}
Let $(X,\calf)$ be a foliated $2$-complex.
  Assume that $\pi_1(X)$ is normally generated by free homotopy classes of curves contained in leaves.

Then $X/\calf$ is an $\bbR$-tree.
\end{thm}

\begin{rem*}
If $f_\sigma$ is constant on $\sigma$, then $\sigma$ is contained in a leaf. This means that contrary to \cite{LP},
we allow 2-simplices to be contained in a leaf. By removing the interior of those $2$-simplices 
(which does not change the leaf space made Hausdorff), 
one can reduce to the case considered by \cite{LP} (for one measured foliation).
For two measured foliations, not allowing them would introduce unnecessary technical complications.
\end{rem*}

\begin{dfn}[Tree dual to a $2$-complex, geometric action]\label{dfn_dual}
Consider a finitely generated group $G$ acting on a tree $T$. We say that $T$ is \emph{dual}
to a foliated $2$-complex $X$ endowed with an action of $G$ if there is an equivariant isometry
between $T$ and $X/\calf$ and if each transverse edge of $X$ isometrically embeds into $X/\calf$.

Then $T$ is \emph{geometric} if it is dual to a foliated $2$-complex $X$ such that the action of 
$G$ on $X$ is free, properly discontinuous, and cocompact.
\end{dfn}

We call a \emph{direct system} of actions on $\bbR$-trees
a sequence of actions
  of finitely generated pairs $G_k\actson T_k$ and an action $G\actson T$,
with epimorphisms $\phi_k:G_k\ra G_{k+1}$ and $\psi_k:G_k\ra G$,
and surjective $\phi_k$-equivariant (resp. $\psi_k$-equivariant) morphisms of $\bbR$-trees 
 $f_k:T_k\onto T_{k+1}$ (resp. $g_k:T_k\onto T$) such that the
  following diagram commutes:
  $$
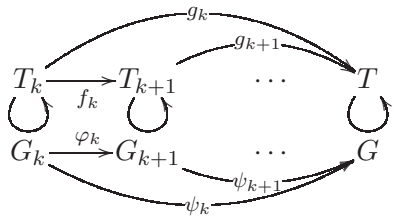
\xymatrix@1@R=0.5cm{T_k\ar[r]_{f_k} \ar@/^0.9cm/[rrr]|{g_k}
    \ar@(dl,dr)[]
    & T_{k+1}  \ar@/^0.5cm/[rr]|{g_{k+1}} \ar@(dl,dr)[] 
    & \cdots
    & T \ar@(dl,dr)[] \\
    G_k\ar[r]^{\phi_k} \ar@/_0.7cm/[rrr]|{\psi_k}
    &G_{k+1} \ar@/_0.4cm/[rr]|{\psi_{k+1}} 
    &\cdots & G}
  $$

For convenience, we will use the notation $f_{kk'}=f_{k'-1}\circ\dots\circ f_k:T_k\ra T_{k'}$.

  \begin{dfn}[Strong convergence]
  We say that a direct system of minimal actions of
  finitely generated groups on $\bbR$-trees $G_k\actson T_k$
  converges strongly to $G\actson T$ if
\begin{itemize*}
\item $G$ is the direct limit of the groups $G_k$
\item for all finite tree $K\subset T_k$, there exists $k'\geq k$ such
  that $g_{k'}$ restricts to an isometry on $f_{kk'}(K)$,
\end{itemize*}
  \end{dfn}

We now recall the definition of a \emph{trivial} strong limit.
If $H$ is a countable group acting by isometries on a metric space $T$, 
denote by $\widehat{T/H}$ the metric space obtained by making Hausdorff the natural pseudo-metric
hold by the quotient space $T/H$.
A strong limit is \emph{trivial} if for $k$ large enough, the space $\widehat{T_k/\ker\psi_k}$ (which is naturally endowed
with an action of $G$) is equivariantly isometric to $T$.

\begin{thm}[{\cite[Corollary 0.3]{LP}}]
  An action of a  finitely generated group $G$ on an $\bbR$-tree $T$ is geometric if and only if every
direct system converging strongly to $T$ converges trivially.
\end{thm}

\subsection{A technical preliminary result}

The following lemma is a slight extension of the result of \cite{LP}
saying that a geometric action is dual to a foliated $2$-complex.

\begin{SauveCompteurs}{lemgeometric}
\begin{lem}\label{lem_geometric}
  Consider a \emph{geometric} action of a finitely generated group $G$ on an $\bbR$-tree $T$, 
and let $X$ be a $2$-complex endowed with a
free properly discontinuous cocompact action of $G$. 
Let $\calf$ be a $G$-invariant measured foliation on $X$.
Consider a map $f:X\ra T$ which in constant on leaves of $\calf$,
and isometric in restriction to transverse edges of $X$.

Then, there exists a $2$-complex $\Hat X$ containing $X$,
 endowed with a free properly discontinuous cocompact action of $G$,
a measured foliation $\Hat\calf$ extending $\calf$, and a map $\Hat f:\Hat X\ra T$ extending $f$,
which is constant on leaves of $\Hat\calf$, and which induces an isometry between
$\Hat X/\Hat\calf$ and $T$.
Moreover, the inclusion $X\subset \Hat X$ induces an epimorphism of fundamental groups.
\end{lem}%
\end{SauveCompteurs}


\begin{proof}
The proof is essentially a rewording of \cite{LP}.
  We choose a large connected finite subgraph $K$ in the $1$-skeleton of $X$, and we describe a construction of a $G$-foliated space $(X_K,\calf_K)$
 containing $X$, and such that the map $X_K/\calf_K\ra T$ is an isometric embedding 
in restriction to the image of $K$. 

First, note that the set $K_0=f(K)$ is a finite tree (\ie the convex hull of finitely many points), 
and has therefore a natural simplicial structure. We can subdivide this simplicial structure so
that for every vertex $v$ of $K$, $f(v)$ is a vertex of $K_0$.
Let $C_K$ be the set obtained by coning off $K$ as follows: glue $K\times [0,1]$ on $K_0$ via
the map $K\times\{1\}\ra K_0$ sending $(x,1)$ to  $f(x)$.
There is a natural measured foliation on $C_K$ induced by the foliation $\{*\}\times [0,1]$ of $K\times[0,1]$.
The set $C_K$ can easily be turned into a simplicial complex whithout subdividing $K\times\{0\}$.
Furthermore, the map $f:X\ra T$ extends uniquely to $f_K:X_K\ra T_K$ as a map constant on leaves
and isometric in restriction to transverse edges.

Let $(X_K,\calf_K)$ be the foliated $2$-complex obtained by gluing on $X$ the set $G\times C_K$
via the map $G\times K\times\{0\}\ra X$ sending $(g,(x,0))$ to $g.x$.
The set $X_K$ has a natural free properly discontinous cocompact action of $G$, so
$X_K$ is a covering of $\ol X_K=X_K/G$. Let $N_K$ be  the image of $\pi_1(X_K)$ in $\pi_1(\ol X_K)$ 
so that $G\simeq \pi_1(\ol X_k)/N_K$.
Note that if $G$ is finitely presented, then one can choose $K$ large enough so that
$\pi_1(K)$ normally generates $\pi_1(X)$, which means in other words that $\Tilde X_K$ is simply connected,
and that $N_K=\{1\}$.
Moreover, by construction, the image of $K$ in $X_K/\calf_K$ isometrically embeds into $T$.

The problem now is that $X_K/\calf_K$ may not be an $\bbR$-tree. In view of Theorem \ref{thm_LP_leaves},
we would need that $\pi_1(X_K)$ is generated by free homotopy classes of curves contained in loops
(note that this is automatically the case for $K$ large enough if $G$ is finitely presented
since one can take $X_K$ to be simply connected).
This is why we are going to consider a Galois covering $\Tilde X_K$ of $\ol X_K$ above $X_K$ so that this condition
is satisfied.

$$\xymatrix@C=0cm@M=0cm{
\Tilde G_K=\pi_1(\ol{X}_K)/\Tilde N_k \ar[d] & \actson &\Tilde X_k \ar[d] \\
G=\pi_1(\ol{X}_K)/ N_k                       &\actson&  X_k \ar[d] \\
                                          &&    \ol{X}_k
}
$$

 Let $\Tilde N_K\subset N_K$ be the normal subgroup of $\pi_1(\ol X_K)$
 generated by free homotopy classes of curves contained in leaves and representing an element of $N_K$.
Let $\Tilde X_K$ be the Galois covering of $\ol X_K$ with deck group $G_K=\pi_1(\ol X_K)/\Tilde N_K$.
Let $\Tilde \calf_K$ be the lift of the measured foliation to $\Tilde X$,
and let $T_K=\Tilde X_K/\Tilde\calf_K$ which is an $\bbR$-tree by Theorem \ref{thm_LP_leaves}.
Denote by $\phi_K:\Tilde G_K\onto G$ the natural epimorphism, 
and let $g_K:T_K\ra T$ be the natural $\phi_k$-equivariant map.
The construction is natural with respect to inclusions $K\subset K'$:
if $K\subset K'$, there is also a natural morphism $\phi_{KK'}:G_K\onto G_{K'}$,
and a $\phi_{KK'}$-equivariant morphism of $\bbR$-trees $g_{KK'}:T_K\ra T_{K'}$.

First, it is clear that $G$ is the direct limit of $G_K$ (with respect to $\phi_{KK'}$ and $\phi_{K}$)
since any relation of $G$ is coned-off in $X_K$ for $K$ large enough.
Moreover, by construction, $T_K$ converges strongly to $T$.
By the defining property of geometric actions, this convergence is trivial.
In other words, 
then for $K$ large enough, the natural map
from $\widehat{T_K/\ker \phi_K}$ to $T$ is an isometry.
But since $\Tilde X_K/\ker \phi_K=X_K$, one has $\widehat{T_K/H_K}=X_K/\calf_K$ so this means that
for $K$ large enough, the map $g_K:X_K/\calf_K\ra T$ is an isometry.
\end{proof}

\subsection{Two foliations on one complex}
Consider a pair of geometric actions $T_1,T_2$.
To exhibit a weak fundamental domain of $\calc(T_1\times T_2)$, our first step
is to write $T_1$ and $T_2$ as the leaf space made Hausdorff of two measured foliations on a common space.

\begin{prop}\label{prop_foliation}
  Let $G$ be a finitely generated group with two geometric actions on $\bbR$-trees
$T_1$, $T_2$. 

Then there is a connected $2$-complex $X$ with a free properly discontinuous cocompact action 
of\/ $G$, and two invariant measured foliations $\calf_1,\calf_2$ on $X$ such that for both $i\in\{1,2\}$,
$T_i$ is dual to $\calf_i$ in the sense of definition \ref{dfn_dual}.

Moreover, if $G$ is finitely presented, we may assume that
$X$ is simply connected.
\end{prop}

Remeber that a map $f$ from a simplex to an $\bbR$-tree $T$ is \emph{affine}
if it is the composition of an affine map to $\bbR$ and of an isometry to $T$.

We will make use of the following standard fact for extending maps to $\bbR$-trees.

\begin{fact}\label{fact_extension}
  Let $X$ be a simplicial $2$-complex, $T$ an $\bbR$-tree, and let $f:X_0\ra T$
be a map defined on a sub-complex $X_0\subset X$ containing the $1$-skeletton of $X$,
and such that $f$ is affine in restriction to each edge.

Then there exists a natural extension $\Hat f:X\ra T$ of $f$,
and a natural subdivision of $X$ which does not change the simplicial structure on $X_0$,
and such that for each simplex $\sigma\subset X_0\setminus X$
(in this new subdivision), $f_{|\sigma}$ is affine. 

If a group $G$ acts on both $X$ and $T$ so that $f$ is equivariant, then the 
subdivision is equivariant, and $\Hat f$ is equivariant. 
\end{fact}

\begin{figure}[htbp]
  \centering
  \input{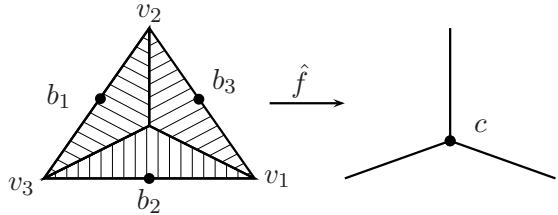}
  \caption{Extension of $f$}
  \label{fig_extension}
\end{figure}

\begin{proof}[Proof of the fact]
  Let $\tau$ be a simplex of $X_0\setminus X$.
We denote by $e_1,e_2,e_3$ the sides of $\tau$,
and by $v_1,v_2,v_3$ its vertices so that $v_i\notin e_i$.
 If $f(\partial\tau)$ contains no tripod (\ie $f(\partial\tau)$
is contained in an arc of $T$), then $f$
can be uniquely extended affinely on $\tau$ (without subdivision).
Otherwise, let $c$ be the center of the tripod $f(\partial\tau)$, and for $i\in\{1,2,3\}$, let $b_i$ its (unique) preimage in $e_i$.
We now subdivide $\tau$ as follows: let $b$ be the barycenter of $\{b_1,b_2,b_3\}$, and
cut $\tau$ along the three segments $[v_i,b]$.
We thus replace the triangle $\tau$ by one new vertex $b$, 3 new edges  $[v_i,b]$, and three new triangles $\{v_i,v_j,b\}$.
We define $\Hat f(b)=c$ and extend $\Hat f$ affinely on the new edges.
But now, $\Hat f$ maps the boundary of the new triangles to intervals, so $\Hat f$ can be uniquely extended
affinely on the new triangles.
\end{proof}

\begin{proof}[Proof of Proposition \ref{prop_foliation}]
 Since $T_1$ is geometric, consider a foliated $2$-complex $X$
to which $T_1$ is dual; more precisely, let $X$
having a free properly discontinuous cocompact action of $G$,
endowed with an invariant measured foliation $\calf_1$, such that $T_1$ is the leaf space 
made Hausdorff of $\calf_1$, and such that transverse edges isometrically embed into $T_1$.
Denote by $f_1:X\ra T_1$ the canonical map.

Note that if $G$ is finitely presented, we can assume that $\pi_1(X)$ is generated by
free homotopy classes of curves contained in leaves (\cite[Remark 2.3]{LP}).
In this case, by gluing finitely many orbits of triangles, and by extending $f_1$ according to
Fact \ref{fact_extension}, one can assume that $X$ is simply connected. 
One has that $T_1$ is dual to the new foliation on $X$ because 
one has two $1$-Lipschitz maps 
$$(X/\calf_1)_{\text{old}}\ra (X/\calf_1)_{\text{new}} \ra T_1$$
such that the composed map $(X/\calf_1)_{\text{old}} \ra T_1$ is an isometry, which forces the two maps to be isometries.

We now define an equivariant map $f:X\ra T_2$.
First, choose an arbitrary equivariant map on the $0$-skeletton (this is possible because the action is free on $X$).
Now extend $f$ affinely on the $1$-skeletton, and use the fact above to extend $f$ on the $2$ skeletton 
(after a subdivision of $X$). Let $\calf_2$ be the measured foliation on $X$ induced by $f$.

Lemma \ref{lem_geometric} shows how to equivariantly enlarge the foliated $2$-complex $(X,\calf_2)$ to a larger $2$-complex 
$(\Hat X,\Hat\calf_2)$ containing $X$, such that the restriction of $\Hat\calf_2$ on $X$ coincides with $\calf_2$,
and such that $T_2$ is dual to $(\Hat X,\Hat\calf_2)$. If $G$ is finitely presented,
since the inclusion $X\subset \Hat X$ induces a epimorphism of fundamental groups, we have that $\Hat X$ is simply connected.

We now extend $f_1$ and $\calf_1$ to $\Hat X$:
first, define $\Hat f_1$ on the set of vertices of $\Hat X\setminus X$ in any equivariant way,
and extend $\Hat f_1$ to the edges of $\Hat X\setminus X$ affinely;
using  fact \ref{fact_extension}, extend $f_1$ and $\Hat\calf_1$ and to the $2$-cells
of $\Hat X\setminus X$ using fact \ref{fact_extension}. 
Since we have natural equivariant $1$-Lipschitz maps  $X/\calf_1\ra \Hat X/\Hat\calf_1\ra T_1$,
and since $X/\calf_1\ra T_1$ is an isometry onto $T_1$, one gets that $ \Hat X/\Hat\calf_1$ is equivariantly isometric
to $T_1$. 
\end{proof}

\subsection{A weak fundamental domain}

We now prove the existence of a weak fundamental domain.

\begin{proof}[Proof of Theorem \ref{thm_fund_domain}]
  Let $X,\calf_1,\calf_2$ be a $2$-complex with two measured foliations as described in Proposition \ref{prop_foliation} above.
Let $f_i:X\ra T_i$ given by the composition of the quotient map $X\ra X/\calf_i$ and the equivariant isometry $ X/\calf_i\ra T_i$.
Let $F=(f_1,f_2):X\ra T_1\times T_2$.
Clearly, $F(X)$ is connected and $G$-invariant. Moreover, the image of each simplex of $X$ is contained in a rectangle.
Since $X$ has finitely many orbits of simplices, $F(X)$ is contained in $G.D$ where $D$ is a finite union of rectangles in $T_1\times T_2$.

Now, the idea is essentially that $F(X)$ should have connected fibers because a fiber should coincide with the image of a leaf.
However, this is not completely true because points of $X/\calf_i$ do not always coincide with leaves of $X$
(we have to \emph{make the space of leaves Hausdorff}). However, by definition of the metric on $X/\calf_i$, 
the following \emph{weak connexity of fibers} holds:
\begin{dfn}[Weak connexity of fibers] Say that a subset $F\subset T_1\times T_2$ has the \emph{weakly connected fibers}
if for all $a,b\in F$ such that $p_i(a)=p_i(b)$, and for each $\eps>0$, there exists a path $c:[0,1]\ra F$ joining $a$ to $b$
in $F$ and such that $p_i\circ c([0,1])$ is contained in the $\eps$-neighbourhood of $p_i(a)$.
\end{dfn}

In view of the characterization of $\calc$ (Proposition \ref{prop_carac_core}, the Theorem follows from the following fact.
\end{proof}

\begin{fact}\label{fact_weak_connexity}
  Let $F\subset T_1\times T_2$ be a connected set enjoying the weak connexity of fibers.
Then its closure $\ol F$ has convex fibers.
\end{fact}

\begin{proof}
  Let $F_V$ be the union of the convex hulls of the vertical fibers of $F$.
Of course, vertical fibers of $F_V$ are convex.
We have $F_V\subset \ol F$. Indeed, let $x\in F_V\setminus F$ so that $x$ lies in a vertical segment
$\{x_1\}\times[a_2,b_2]$ whose enpoints are in $F$. Let $\eps>0$ and $c$ be a path in $F$ joining 
$(x_1,a_2)$ to $(x_1,b_2)$, and such that $p_1\circ c$ is contained in the $\eps$-neighbourhood of $\{x_1\}$.
Then $p_2\circ c([0,1])$ must contain $x_2$. 
This implies that $x$ lies in the closure of its
horizontal fiber $F\cap\, T_1\times\{x_2\}$, and in particular that $F_V\subset\ol F$.
This also implies that connected horizontal fibers of $F$ stay connected in $F_V$, 
that $F_V$ is connected and still satisfies the weak connexity of (horizontal) fibers.

Now let $F_{V,H}$ be the union of the convex hulls of the horizontal fibers of $F_V$.
The symmetric argument shows that horizontal and vertical fibers of $F_{V,H}$ are convex,
and that $F_{V,H}\subset\ol{F}$.
Corollary \ref{cor_adherence} concludes that $\ol F$ has convex fibers.
\end{proof}

\subsection{Finiteness of intersection number}

We now prove the finiteness of the intersection number of two geometric actions of a finitely presented group.

\begin{proof}[Proof of Proposition \ref{prop_finite_intersection}]
 Let $X,\calf_1,\calf_2$ be a $2$-complex with two measured foliations 
 such that  $X/\calf_i$ is equivariantly isometric to $X/\calf_i\ra T_i$
as described in Proposition \ref{prop_foliation}.
Let $f_i:X\ra T_i$ given by the composition of the quotient map $X\ra X/\calf_i$ and the equivariant isometry $ X/\calf_i\ra T_i$.
Let $F=(f_1,f_2):X\ra T_1\times T_2$.

Since $G$ is finitely presented, we can assume that $X$ is simply connected.
Therefore, we can apply Lemma 3.4 of \cite{LP} which claims that the natural maps from the set of leaves of $\calf_i$
to $X/\calf_i$ is one-to-one outside a countable set.
In particular, all but countably many fibers of $F(X)$ are connected.
Finally, $F(X)$ is connected,
$G$-invariant and contained $G.D$ where $D$ is a finite union of rectangles in $T_1\times T_2$.
The proposition then follows from Lemma \ref{lem_countable}.
\end{proof}

\begin{lem}\label{lem_countable}
  Let $F\subset T_1\times T_2$ be a connected set, enjoying the weak connexity of fibers. 
Assume that all but countably many horizontal and vertical fibers of $F$ are convex.

Then $\ol F$ have convex fibers and $\ol F\setminus F$ has measure $0$.
\end{lem}

\begin{proof}
The first assertion follows from Fact \ref{fact_weak_connexity}.
  Let $F_V$ be the union of the convex hulls of the vertical fibers of $F$.
Since $F_V\setminus F$ is contained in countably many fibers, 
it has null measure.
Moreover, the proof of Fact \ref{fact_weak_connexity} shows that connected horizontal
fibers of $F$ stay connected in $F_V$, so all but countably many $F_V$ horizontal
fibers of $F_V$ are connected. Now let $F_{V,H}$ be the union of the convex hulls of the horizontal fibers of $F_V$.
The symmetric argument shows that horizontal and vertical fibers of $F_{V,H}$ are convex.
Since $F_{V,H}\setminus F$ has null measure, the lemma is a consequence of the following lemma.
\end{proof}

\begin{lem}
    Let $F\subset T_1\times T_2$ be a connected set with convex fibers.
Then $\ol F\setminus F$ has measure $0$.
\end{lem}

\begin{proof}
  We have to prove that the trace of $\ol{F}\setminus F$ on each compact rectangle has measure zero.
Let $R_0$ be a compact rectangle, and let $J_i=p_i(\ol{F}\cap R_0)$.
Let $R=J_1\times J_2$. Since $\ol{F}\cap R_0\subset R$, we have to prove that 
$(\ol{F}\setminus F)\cap R$ has measure 0.
By Lemma \ref{lem_fibre_quadrant}, $\ol{F}$ is the complement of a union of quadrants.
By changing $T_1$ and $T_2$ to $p_1(\ol{F})$ and $p_2(\ol{F})$, we may also assume that
$p_i(\ol{F})=T_i$. Thus, one can apply lemma \ref{lem_contractile_general} saying that the semi-flow $\phi_t$ 
contracting $T_1\times T_2$ onto $\ol{F}$
is well defined and that $R$ is $\phi_t$-invariant. 

We now prove that $\ol{F}\setminus F$ is contained in $\phi_\infty(\partial R)$.
Since $\phi_\infty$ is a Lipschitz map, this will conclude the proof of the fact.
Let $x=(x_1,x_2)\in (\ol F\setminus F)\cap R$. We can assume that $x\notin \partial R$.
Consider the four open rectangles obtained from $R$ by subdividing $R$ vertically and horizontally at $x$.
At least one of those open rectangles, call it $S$, does not intersect $\ol F$, since otherwise, $F$ would
intersect the corresponding four quadrants, which would imply $x\in F$ by Fact \ref{fact_quadrants}.
Starting from $x$, one can follow in $S$ the semi-flow $\phi_t$  backwards until we reach the boundary of $R$
at a point $x_0$, so that $\phi_{t_0}(x_0)=x$ for some $t_0\in\bbR_+$.
Now $\phi_\infty(x_0)=\phi_\infty(x)=x$ since $x\in \ol{F}$.
\end{proof}

\section{Application to automorphisms of free groups}\label{sec_automorphisms}

In this section, we study in more detail the case of attracting
and repulsive trees of an outer automorphism of a free group $G$.

Our starting point is the following theorem due to Levitt and Lustig:

\begin{thm}[Levitt-Lustig]
  Let $\alpha$ be an automorphism of a finitely generated free group $G$.
Then there exists two actions of $G$ on $\bbR$-trees $T_1,T_2$ and two 
homotheties $h_1:T_1\ra T_1$, and $h_2:T_2\ra T_2$ with stretching fators $\lambda_1,\lambda_2$ 
such that
\begin{itemize*}
\item $h_i$ is $\alpha$-equivariant (\ie $h_i(g.x)=\alpha(g).h_i(x)$)
\item $\lambda_1\geq 1$, $\lambda_2\leq 1$
\item the action of $G$ on $T_i$ has trivial arc fixators 
\item $g\in G$ is elliptic in $T_1$ if and only if it is elliptic in $T_2$
\item There exists $\eta>0$ such that if $l_{T_1}(g)+l_{T_2}(g)\leq \eta$, then $g$ is elliptic in $T_1$ and $T_2$
(so  $l_{T_1}(g)+l_{T_2}(g)=0$).
\end{itemize*}
\end{thm}

\newcommand{\iwip}{iwip}
Note that this implies that the action of $G$ on $T_1\times T_2$ is discrete.
Say that an automorphism $\alpha$ is irreducible with irreducible powers (\emph{\iwip})
no free factor of $G$ is periodic under $\alpha$ up to conjugacy.
The theorem above was proved for \iwip\ automorphisms in \cite{Lus_discrete} and \cite{BFH_laminations,BFH_erratum_laminations}.
Moreover, in this case, either the action of $G$ on $T_1$ and $T_2$ is free,
or $\alpha$ is induced by a pseudo-ansov homeomorphism of a surface with boundary.

The goal of this section is the following theorem.

\begin{thm}\label{thm_automorphism}
  Let $\alpha$  be an automorphism of a finitely generated free group $G$, and let $T_1,T_2$ be as above.
Assume that the actions of $G$ on $T_1$ and $T_2$ are geometric, and that $\lambda_1>1$.
Then there is a $G$-invariant, cofinite (discrete) set $\cals\subset\calc$ such that
\begin{itemize*}
\item $\calc$ has a structure of simplicial complex, and $\calc\setminus \cals$ is a surface
\item the link in $\calc$ of each point of $\cals$ is a disjoint union of lines and circles
\item the action of $G$ on $\calc$ is discrete and cocompact, and it is properly discontinuous on $\calc\setminus\cals$
\item $\Sigma=\calc/G$ is a \emph{pinched} compact surface, and $\Sigma\setminus(\cals/G)$
is a surface with finitely many puntures
\item $\Sigma$ has two transverse measured foliations $\calf_1,\calf_2$,
and the map $H=(h_1,h_2)$ induces a homeomorphism of $\Sigma$ preserving $\cals/G$,
preserving both foliations, and which multiplies the transverse measure of $\calf_i$ by $\lambda_i$,
and $\lambda_2=1/\lambda_1$.
\end{itemize*}
\end{thm}

In this theorem, a \emph{pinched surface} is a space obtained from a compact surface
by collapsing to a point finitely many finite sets. Equivalently, this is a finite $2$-complex
in which the link of every vertex is a finite disjoint union of circles.
\begin{SauveCompteurs}{cpt_iwip}
\begin{cor}\label{cor_iwip}
  Assume that $\alpha$ is irreducible with irreducible powers.
If $T_1$ and $T_2$ are both geometric, then $\alpha$ is induced by
a pseudo-anosov homeomorphism of a surface with boundary.
\end{cor}
\end{SauveCompteurs}

\subsection{A structure of 2-complex}

  By Theorem \ref{thm_fund_domain}, there is a weak fundamental domain for $\calc$:
there is a finite set of rectangles $\cald$ such that, the closure of the set $X=\bigcup_{R\in\cald,g\in G} g.R$
is connected and has convex fibers, which implies that $\calc\subset \ol{X}$.

Our first goal is to prove that $X$ has a $G$-invariant structure of $2$-complex,
and that the action of $G$ on $X$ is properly discontinuous outside $\cals_0$ (Corollary \ref{cor_prop_discont}).

First, $\calc$ is non-empty since $T_1$ and $T_2$ are necessarily irreducible.
Moreover, there is no twice light rectangle: indeed, since there is no arc in $T_1$ containing no branch point
since the image of such an arc under $h_1$ would give arbitrarily large such arcs, which is not possible
since by minimality, there is a finite tree whose translates under $G$ cover $T_1$.\\

On $T_1\times T_2$, we use the distance $d((x_1,x_2),(y_1,y_2))=d_{T_1}(x_1,y_1)+d_{T_2}(x_2,y_2)$.
We also denote $l(g)=l_{T_1}(g)+l_{T_2}(g)$.
Let $\eta$ such that  $l(g)<\eta \Rightarrow l(g)=0$.
Denote by $\cals_0\subset T_1\times T_2$ the set of points
having a non-trivial stabilizer under the action of $G$.

\begin{lem}
  Any two distinct points $x,x'\in \cals_0$ are at distance at least $\eta/2$.
\end{lem}

\begin{proof}
  let $g,g'\in G\setminus\{1\}$ fixing $x$ and $x'$ respectively.
Since arc fixators are trivial, $\Fix g$ and $\Fix g'$ are reduced to a point in both $T_1$ and $T_2$.
But the translation length of $gg'$ in $T_i$ is $2d_{T_i}(\Fix g,\Fix g')$.
Therefore, $l(gg')=2d(x,x')$, and by hypothesis, we get $d(x,x')=0$ or $d(x,x')>\eta/2$.
\end{proof}

\begin{lem}
  If $d(x,g.x)<\eta$ then $g$ has a fix point $y$ such that $d(x,y)= d(x,g.x)/2$.
\end{lem}

\begin{proof}
  Since $d(x,g.x)<\eta$, $l(g)<\eta$ so $g$ is elliptic in $T_1$ and $T_2$.
This implies that $g$ fixes the midpoint of $[x_i,g.x_i]$ in $T_i$.
\end{proof}

\begin{lem}
The set  $X$ is closed in $T_1\times T_2$.
\end{lem}

\begin{proof}
Let $x=(x_1,x_2)\in \ol{X }$, and consider a sequence $x_i$ in $\bigcup_{R\in \cald} R$
and a sequence $g_i\in G$ such that $g_i.x_i$ converges to $x$. Up to extracting a subsequence, we can
also assume that $x_i$ converges to $a=(a_1,a_2)$. This implies that $g_i.a$ converges to $x$.
If for some $i$ and infinitely many $j>i$, $g_jg_i\m=1$, then $x\in X$, and we are done.
Otherwise, given $\eps>0$, for $i$ large enough and all $j>i$,
$g_jg_i\m$ moves the point $g_i.a$, and thus $x$ by less that $\eps$.
Therefore, all the elements $g_jg_i\m$  have a common fix point $b=(b_1,b_2)\in \cals_0$ with $d(b,x)\leq \eps$.
Since this is valid for all $\eps > 0$, one gets $b=x$. Therefore, $d(g_j.a,x)$ is constant so it must be 0, and $x\in X$.
\end{proof}

From now on, we subdivide $\cald$ so that each rectangle $R\in\cald$ has diameter at most $\eta/10$,
and so that $R\cap\cals_0$ is either empty or consists of a single corner of $R$.

\begin{lem}\label{lem_proper}
  $X\setminus \cals_0$ is locally compact, and $G$ acts properly discontinuously on $X\setminus \cals_0$.

More precisely, every point of $X\setminus \cals_0$ has a neighbourhood $V$ such that for each $R\in \cald$,
there is at most one $g\in G$ such that $g.R$ intersects $V$.
\end{lem}

\begin{proof}
Clearly, the second part of the statement implies the fact that $X\setminus\cals_0$ is locally compact and that
the action of $G$ on $X\setminus\cals_0$ is properly discontinuous.

Let $x\notin\cals_0$ and let $V=V_1\times V_2$ be a product of balls around $x$ 
of radius at most $\eta/10$ and such that $V$ does not intersect $\cals_0$.
Fix a rectangle $R\in \cald$ and assume that $g_1.R$, $g_2.R$ intersect $V$ for some $g_1\neq g_2\in G$.
Then $h=g_2g_1\m$ moves a point of $g_2.R\cap V$ by at most $\eta/2$, so $h$ has a fix point $y=(y_1,y_2)$ at distance at most
$\eta/2$ of $x$. Since $R$ meets $\cals_0$ only at corners of $R$, $g_1.R$ is contained in 
a quadrant $\delta_1\times\delta_2$  based at $y$. Now since $y\notin V$, one may assume that $y_1\notin V_1$,
so $V_1$ is contained in $\delta_1$.
Since the action of $\Stab\{y_1\}$ on the set of directions based at $y_1$ is free, 
$p_1(h.R)\subset h.\delta_1$ implies that $h.R$ does not intersect $V$.
\end{proof}

\begin{lem}\label{lem_2complex}
There is a subdivision of $\cald$ and a $G$-invariant cofinite set $\cals$ such that
\begin{itemize*}
\item two rectangles of $G.\cald$ either coincide or intersect in a (maybe degenerate or empty) interval, 
contained in a horizontal or vertical side of $R$,
and whose endpoints are in $\cals$;
\item for each rectangle $R\in\cald$,
$\rond{R}$ is open in $X$
\end{itemize*}

Moreover, $X$ has a $G$-equivariant structure of a simplicial complex with finitely many orbits of cells,
and whose vertex set is $\cals$.
\end{lem}

\begin{rem*}
  The rectangles of $\cald$ don't make $X$ a rectangle complex. In general, there is no decomposition
of $X$ as a rectangle complex whose sides are vertical and horizontal, as can be seen in
the torus endowed with two irrational foliations.
\end{rem*}

\begin{proof}
  For each rectangle $R=I_1\times I_2\in \cald$, let $\cali(R)$ be the (finite) set of 
translates of rectangles $R'\in \cald$ which intersect $R$ in a non-degenerate rectangle.
Let $\cali_i(R)\subset I_i$ 
be the finite set of endpoints of the projections on $I_i$ of the rectangles $R\cap R'$ for $R'\in\cali(R)$.
After
subdividing each rectangle $R=I_1\times I_2\in \cald$ along $\cali_1(R)\times I_2 \cup I_1\times \cali_1(R)$,
one gets that $g.R'\cap R$ is a segment or a point whenever $g.R'\neq R$.

Let $\cals$ be the set of points of $X$ occuring as an endpoint (or a corner) of 
the intersection of two elements of $G.\cald$.
Clearly, $\cals_0\subset\cals$, and $\cals$ intersects each rectangle of $\cald$ in a finite set so 
 $\cals/G$ is finite.

The fact that $\rond R$ is open in $X$ follows from the fact that two rectangles intersect only in their boundary
and that a point of $\rond R$ has a neighbourhood intersecting only finitely many rectangles.

Finally, we triangulate $X$ as follows: given $R\in\cald$,
choose a triangulation of $R$ 
whose vertex set is $R\cap\cals$ (this contains the corners of $R$). 
Since given $R,R'\in \cald$ there is at most one element $g\in G$ 
such that $g.R=R'$, this triangulation can be chosen so that it extends $G$-equivariantly to $X$.
\end{proof}

\subsection{The core as a subcomplex}

In this section, we prove that up to subdividing our structure of $2$-complex on $X$,
$\calc$ appears as a subcomplex of $X$.

Roughly speaking, the following lemma says that $X$ 
cannot go in a direction branching from the interior of a rectangle.

\begin{lem}
  Let $R\in \cald$, let $x\in\rond R$, and let $Q$ be a quadrant based at $x$, disjoint from $R$.
Then $Q$ is disjoint from $X$.
\end{lem}

\begin{proof}
Assume on the contrary that there is a point $y=(y_1,y_2)\in X\cap Q$.
Write $Q=\delta_1\times\delta_2$, and let $(b_1,b_2)$ be the base point of $Q$.
Choose a path $c=(c_1,c_2):[0,1]\ra X$ joining $y$ to $x$ in $X$, and let $t_0$ be the first time for which
one coordinate of $c(t_0)$ equals $b_1$ or $b_2$. Note that for $t<t_0$, $c(t)\notin R$.
Since $\rond R$ is open in $X$, $c(t_0)\neq (b_1,b_2)$.
So assume for instance that $c_1(t_0)=b_1$ and  $c_2(t_0)\in \eta_2$.
By convexity of fibers of $X$, $\{b_1\}\times [b_2,c_2(t_0)]\subset X$,
which contradicts the fact that $\rond R$ is open in $X$.
\end{proof}

We now study how $\calc$ may intersect rectangles of $\cald$.

\begin{lem}
One can subdivide equivariantly each rectangle of $\cald$ into smaller subrectangles such that for all $R\in \cald$,
either $R\subset \calc$ or $R\cap\calc\subset \partial R$.
\end{lem}

\begin{proof}
Let $R\in \cald$ be a rectangle which is not contained in $\calc$.
We are going to prove that $\calc\cap R$ can be obtained from $R$ by removing a finite number of quadrants.
The lemma will follow easily after subdivision of $R$.

Let $Q$ be a light quadrant. Then one of the following holds:
\begin{enumerate*}
\item\label{type_nul} either $\rond R\subset Q$, or $R\cap Q=\es$
\item\label{type_demi} or $R\cap Q$ is an open half rectangle
\item\label{type_coin} or $R\cap Q$ contains exactly one corner of $R$
\end{enumerate*}

We are not interested in quadrants of type \ref{type_nul}.
We also can assume that up to subdividing $R$, we can avoid quadrants of type \ref{type_demi}.
Indeed, let $I$ be a side of $R$, and let $L_{I}\subset R$ be the union of the traces on $R$ of all light quadrants of type \ref{type_demi}
containing $I$. If $L_{I}$ contains $\rond R$, then $\calc$ does not intersect $\rond{R}$ and we are done.
Otherwise, we can cut $R$ into $\ol L_I$ and $R\setminus L_I$. Doing this for the four sides of $R$,
we can avoid the occurence of quadrants of type \ref{type_demi}.

Note that a quadrant of type \ref{type_coin} is necessarily based at a point $x\in R$.
Now let $c=(c_1,c_2)$ be a corner of $R$, and let $L_{c}$ be the unions of the traces on $R$ 
of all light quadrants of type \ref{type_coin} containing $c$,
and assume that $L_c\neq \es$.
We are going to prove that $L_{c}$ is itself the trace of a light quadrant of type \ref{type_coin}.
The lemma will follow since $R\cap\calc$ will be obtained from $\calc$ by removing at most 4 quadrants.

Denote by $]a,c]=]a_1,c_2]\times\{c_2\}$ and $]b,c]=\{c_1\}\times ]b_2,c_2]$ the intersection of $L_{c}$ with
the two sides of $R$ containing $c$.
Let $P=\rho_1\times\rho_2$ be the quadrant based at $(a_1,b_2)$ and containing $c$.
Clearly, any light quadrant $Q=\delta_1\times \delta_2$ of type \ref{type_coin} is contained in $P$:
for instance, $a\notin Q$ and $c\in Q$ implies that $\delta_1\subset\rho_1$.
Thus $L_c\subset P$.

We now prove the other inclusion. Let $x=(x_1,x_2)\in P$. 
Since $x_1\in ]a_1,c_1]$, by definition of $a$, there exists a light quadrant $Q=\delta_1\times\delta_2$
with $x_1\in\delta_1$.
Similarly, there exists  a light quadrant $Q'=\delta'_1\times\delta'_2$
with $x_2\in\delta'_2$. 
It suffices to prove that the quadrant $Q''=\delta_1\times\delta'_2$ is light.
We choose a base point $*\in X$, and consider a sequence $g_k.*\in Q''$ making $Q''$ heavy.
Then for $k$ large enough, $g_k.*\notin Q$, $g_k.*\notin Q'$ and $g_k.*\notin R$.
It follows that $g_k.*$ lies in a quadrant branching from $R$ as in the previous lemma.
Since $g_k.*\in X$, this is a contradiction.
\end{proof}

Let $\cald'$ be the set of rectangles $R\in\cald$ which are contained in $\cald$
together with the set of horizontal or vertical intervals obtained
as the trace of $\calc$ on the boundary of those rectangles (we think of those segments
as degenerate rectangles).

Using a triangulation as in Lemma \ref{lem_2complex}, we immediately get the following corollary:

\begin{cor}\label{cor_prop_discont}
  The core $\calc$ of $T_1\times T_2$ has a $G$-invariant structure of $2$-complex,
$\calc\setminus\cals_0$ is locally compact, and the action of $G$ on $\calc\setminus \cals_0$
is properly discontinuous.
\end{cor}

\subsection{Singular points of the core}

We now study  the $2$-complex $\Sigma=\calc/G$, and prove that it's a pinched surface.
This will conclude the proof of Theorem \ref{thm_automorphism}.
To this means, we are going use the map $H=(h_1,h_2):T_1\times T_2\ra T_1\times T_2$.

\begin{lem}
  The core $\calc$ is $H$-invariant. Moreover, $H$ induces an homeomorphism $\ol H$ of $\Sigma$,
and $\lambda_1\lambda_2=1$.
\end{lem}

Note that two projections $p_i:\calc\ra T_i$ induce two natural measured foliations
on $\calc$ and $\Sigma$, which are preserved by $H$ and $\ol H$ respectively. Moreover, $\ol H$ expands
the  corresponding transverse measures by a factor $\lambda_2$ and $\lambda_1$ respectively.

\begin{proof}
The first part follows from the fact that $H$ sends quadrant to quadrant,
or from the characterization of $\calc$.
Moreover, $x$ and $y$ are in the same $G$-orbit if and only if $H(x)$ and $H(y)$ are.
The transverse measures of the two foliations on $\Sigma$ define a finite non-zero
measure on $\Sigma$: if the measure was 0, then $\calc$ would be a simplicial tree, implying that $T_1$ and $T_2$ are
both simplicial, which is impossible since $\lambda_1\neq 1$. 
Since $\ol H$ expands this measure by a factor $\lambda_1\lambda_2$,
it follows that $\lambda_1\lambda_2=1$.
\end{proof}

Let $\cals$ be the set of points of $\calc$ occuring as an endpoint (or a corner) of 
the intersection of two elements of $G.\cald'$.
Let $\ol\cals$ be its image in $\Sigma$.

\begin{lem}
$\calc\setminus \cals$ is a surface,
and the link of each point of $\cals$ in $\calc$ is a disjoint union of circles and lines. 
Moreover, $\Sigma\setminus\ol\cals$ is a closed surface with finitely
many punctures, and $\Sigma$ is a pinched surface.
\end{lem}

\begin{proof}
Each point in $\calc\setminus\cals$ has a neighbourhood of the form $I_1\times K_2$
or $K_1\times I_2$ where $K_i\subset T_i$ is a finite subtree and $I_i\subset T_i$ is a compact segment.
Let $\cals_H$ (resp.\ $\cals_V$) be the \emph{horizontal} (resp.\ \emph{vertical}) \emph{singular set} defined as
the set of points $x\in\calc\setminus\cals$ having 
a neighbourhood of the form $K_1\times I_2$ (resp. $I_1\times K_2$) but no neighbourhood of the form $I_1\times I_2$.

 We aim to prove that $\cals_H$ and $\cals_V$ are empty.
The image $\ol\cals_H$ of $\cals_H$ in $\Sigma$ is a finite union of horizontal open edges,
and has therefore a finite horizontal transverse measure, which is non-zero if $\cals_H\neq\es$.
Now the homeomorphism $\ol H$ preserves $\cals_H$ and expands its horizontal transverse measure by $\lambda_1>1$.
This only allows $\cals_H=\es$. Similarly, $\cals_V=\es$, so $\calc\setminus \cals$ is a surface.  

Since the action of $G$ on $\calc\setminus\cals$ is free and properly discontinuous, 
this implies that $\Sigma\setminus\ol\cals$ is a surface.
Therefore, the link at any point in $\ol\cals$ is a graph having only valence $2$-vertices.
Since it is compact, it is a finite disjoint union of circles.
It follows that $\Sigma\setminus\ol\calc$ is a closed surface with finitely many punctures,
and that $\Sigma$ is a pinched surface.
Similarly, the link in $\calc$ of any point of $\cals$ is a disjoint union of lines and circles.
\end{proof}

Theorem \ref{thm_automorphism} follows.

\subsection{The case of an irreducible automorphism with irreducible powers}

We prove here Corollary \ref{cor_iwip}, \ie that if $\alpha$
does not preserve any free factor of $G$ up to conjugacy,
and if $T_1$ and $T_2$ are geometric, then $\alpha$ is induced by a pseudo-Anosov
automorphism on a surface with boundary.

\begin{proof}[Proof of Corollary \ref{cor_iwip}]
We make use the fact that if $\alpha$ is not induced by a 
pseudo-Anosov automorphism on a surface with boundary, then the actions of $G$ on $T_1$ and $T_2$ are free (\cite{BH_tt}).
We assume that $T_1$ and $T_2$ are geometric and argue towards a contradiction.
Since the action on $T_1$ and $T_2$ is free, $\calc$ is the universal cover of $\Sigma$, and $G\simeq\pi_1(\Sigma)$.

Up to removing some points from $\ol\cals$, we can assume that the link of each point of $\ol\cals$ is not a circle.
If $\cals=\es$ we have a contradiction as $G$ is a not a surface group.
Let $\Hat\Sigma$ be the $2$-complex obtained by blowing-up $\Sigma$ as follows:
for each vertex $v\in\ol\cals$, consider the tree $T_v$ defined as the 
complete bipartite graph on $\{v\}$
and $C(v)$ where $C(v)$ is the set of connected components of the link of $v$;
then glue back $C(v)$ to $\sigma\setminus\ol\cals$.
What we get is a (maybe not connected) surface together with a finite collection of finite trees
attached to finitely many points. Note that $\Hat\Sigma$ is homotopy equivalent to $\Sigma$.
However, since $G$ contains no non-trivial surface group ($G$ is free),
all those surfaces are spheres. But this contradicts the fact that the universal covering $\calc$ of $\Sigma$
is contractible.
\end{proof}

\section{Equality with Scott and Swarup's intersection number}\label{sec_scott}

Let $T_1$ and $T_2$ be two actions of a finitely generated group $G$ on simplicial trees with one orbit of edges.
Scott defined in \cite{Sco_symmetry} the intersection number of $T_1$ and $T_2$.
One can reword the definition as follows. Choose $e_1,e_2$ two oriented edges in $T_1$ and $T_2$.
Let $\delta(e_i)\subset T_i$ be the direction based at the origin of $e_i$, and containing $e_i$.
Let $H_i$ be the stabilizer of $e_i$ which is also the stabilizer of $\delta(e_i)$.
Fix a base point $*_i\in T_i$.
Let $X_{e_i}=\{g\in G| g.*_i\in\delta(e_i)\}$. $X_{e_i}$ is clearly left invariant under 
$H_i$.
Let $C$ be a Cayley graph of $G$.
In Scott's terminology, $X_{e_i}$ is a $H_i$-almost invariant set, 
which means that $X_{e_i}/H_i$ has a finite coboundary in the graph $C/H_i$.
Note that changing the basepoint above only changes $X_{e_i}/H_i$ by a finite set, and such
a change will be not be important in the following definitions.

Scott says that $X_{e_1}$ and $X_{e_2}$ \emph{cross} if the four sets
$X_{e_1}\cap X_{e_2}$, $X_{e_1}^*\cap X_{e_2}$, $X_{e_1}\cap X_{e_2}^*$, $X_{e_1}^*\cap X_{e_2}^*$ project to infinite
sets in $C/H_1$.
They prove that this occurs if and only if those four sets also
project to infinite sets in $C/H_2$.
Note that the fact that $X_{e_1}$ crosses $X_{e_2}$ does not change if we replace one edge with the edge with
the reverse orientation. Moreover, $X_{e_1}$ crosses $X_{e_2}$ if and only if $X_{g.e_1}$ crosses $X_{g.e_2}$.

Finally, $i(T_1,T_2)$ is defined as the number of double cosets $H_1 g H_2$
such that $X_{e_1}$ crosses $g.X_{e_2}$. 
Now the set of double cosets is in one-to-one correspondance with the set
of orbits of pairs of non-oriented edges in $T_1\times T_2$.
In other words, $i(T_1,T_2)$ is the number of orbits pairs of non-oriented
edges $(e_1,e_2)$ in $T_1\times T_2$ such that $X_{e_1}$ crosses $X_{e_2}$.

\begin{prop}
  Our definition of the intersection number coincides with Scott's definition.
\end{prop}

\begin{proof}
To identify our definition of intersection number with Scott's one,
we just need to check that a rectangle $e_1\times e_2$ is contained in $\calc$ if and only if $X_{e_1}$ and $X_{e_2}$ cross.
Denote by $Q(e_1,e_2)$ the quadrant $\delta(e_1)\times\delta(e_2)$.

Assume first that $X_{e_1}$ does not cross $X_{e_2}$, and let's prove that $e_1\times e_2$ is not contained in $\calc$. 
Up to a good choice of orientations of $e_1$ and $e_2$,
one gets that the image of $X_{e_1}\cap X_{e_2}$ projects to a finite set in $G/H_1$.
This means that the set
$Z=\{g.*_1 | g.(*_1,*_2)\in Q(e_1,e_2) \}$ meets only finitely many $H_1$-orbits.
In particular, $Z$ is bounded, so $Q(e_1,e_2)$ is light.
Since the open rectangle $e_1\times e_2$ is contained in $Q(e_1,e_2)$,
we get that $e_1\times e_2$ is not contained in $\calc$.

For the converse, assume that $X_{e_1}$ cross $X_{e_2}$,
and let's check that the rectangle $e_1\times e_2$ lies in $\calc$.
We first treat the case where $\calc$ is non-empty, and we choose a base point $*\in\calc$.
Assume by contradiction that there is a light quadrant $Q$ containing $e_1\times e_2$.
In this case, there is a choice of orientations of $e_1$ and $e_2$ such that
$Q$ contains $Q(e_1,e_2)$. In particular, $Q(e_1,e_2)$ is light. 
But since $X_{e_1}$ crosses $X_{e_2}$ there is an element $g\in G$ sending
$*$ in $Q(e_1,e_2)$. Since $Q(e_1,e_2)$ is light, this contradicts the fact that
$g.*\in\calc$.

We now check that if $\calc$ is empty, then Scott's intersection number is zero.
We want to prove that for each pair of non-oriented edges $(e_1,e_2)$, $X_{e_1}$
does not cross $X_{e_2}$. This is easy for homothetic actions on a line.
Otherwise, Remark \ref{rem_coeurvide} shows that there is a good choice of orientations
for $(e_1,e_2)$, and a good choice of base point $*$ such that the orbit of $*$
does not meet $Q(e_1,e_2)$; this implies that $X_{e_1}$ does not cross $X_{e_2}$.
\end{proof}

\section{Strong intersection number and asymmetric core}\label{sec_strong_asymmetric}

\subsection{Asymmetric core}
We first introduce an \emph{aymmetric core} (actually two asymmetric cores $\cala_1,\cala_2$) 
for two actions of $G$ on simplicial trees.
Note that there is no obvious generalisation of $\cala_1$ when $T_1$ is an $\bbR$-tree.

In what follows, given an action of $G$ on a tree $T$, we denote by $G_x$ the stabilizer of a point or an edge $x$ of $T$.

\begin{dfn}[Asymmetric core]
  Assume that a group $G$ acts on two simplicial trees $T_1,T_2$. Assume that each edge stabilizer of $T_1$
acts non-trivially on $T_2$.

Then the left asymmetric core $\cala_1(T_1\times T_2)$ is defined by
$$\cala_1=\bigcup_{x_1\in T_1}\{x_1\}\times\Min_{T_2}(G_{x_1}).$$

The right asymmetric core $\cala_2(T_1\times T_2)$ is defined symmetrically under symmetrical hypotheses.
\end{dfn}

\begin{rem}
It is clear that $\cala_1$ is the smallest invariant subset of $T_1\times T_2$ which has non-empty connected $p_1$-fibers.
The fact that each edge stabilizer of $T_1$ acts non-trivially on $T_2$ implies that 
$\calc$ is non-empty and that $T_1$ and $T_2$ are not refinements of a common simplicial action.
In particular, $\cala_1$ and $\cala_2$ are contained in $\calc$ by Proposition \ref{cor_carac_core}.
Moreover, $\cala_1$ is contractible since it maps to a tree with contractible fibers.
\end{rem}

More generally, if one only assume that each edge or vertex stabilizer of $T_1$ either acts non-trivially, or has a fix point
in $T_2$, then one can do a similar, but less canonical construction:
for each vertex or edge $x_1\in T_1$ such that $G_{x_1}$ acts non-trivially on $T_1$, we still 
impose that $\cala_1\cap p_1\m(x_1)=\{x_1\}\times \Min_{T_2}(G_{x_1})$;
for ${x_1}\in V(T_1)$ such that $G_{x_1}$ fixes a point in $T_2$, we choose a fix point $x_2$ of $G_{x_1}$ in $T_2$, and define
$\cala_1$ so that $\cala_1\cap p_1\m(x_1)=\{(x_1,x_2)\}$;
For each edge $e_1\in E(T_1)$ with endpoints $a_1,b_1$, and such that $G_{e_1}$ is elliptic in $T_2$, 
we choose a segment $I$ in $p_1\m(e_1)$ joining a fix point of $G_{e_1}$ in $p_1\m(a_1)$
to  a fix point of $G_{e_1}$ in $p_1\m(b_1)$ (since this segment might be oblique, $\cala_1$ is not in general a subcomplex of $T_1\times T_2$).
In this more general situation, it remains true that $\cala_1$ is contractible, but it depends on choices. Note however that
the $2$-cells of $\cala_1$ are independant of choices since they come from the minimal subtrees (in $T_2$) of edge stabilizers (in $T_1$).

\subsection{Strong intersection number.}

Given two actions of a finitely generated group $G$ on two simplicial trees $T_1,T_2$ with one orbit of edges,
Scott and Swarup give a definition of the strong intersection number, which we can reword as follows.
Recall that given an oriented edge $e_i\in T_i$, $\delta(e_i)$ denotes
 the direction based at the origin of $e_i$ and containing $e_i$, and that
$X_{e_i}=\{g\in G| g.*_i\in\delta(e_i)\}$. 
Let $\partial X_{e_i}$ be the boundary of $X_{e_i}$ in a Cayley graph $C$ of $G$ with respect to some finite generating set $S$,
\ie $$X_{e_i}=\{g\in X_{e_i}|\exists s\in S\cup S\m\, gs\notin X_{e_i}\}.$$

Scott and Swarup say that $X_{e_1}$ \emph{crosses strongly} $X_{e_2}$  if both
$\partial X_{e_1}\cap X_{e_2}$ and $\partial X_{e_1}\cap X_{e_2}^*$ project to infinite
sets in $C/G_{e_2}$. This notion does not depend on the choice of $S$.

\begin{lem}
Consider $T_1,T_2$ two simplicial trees dual to one edge splittings a finitely generated group $G$.
Let $e_1,e_2$ be two edges in $T_1$ and $T_2$ respectively.

Then $X_{e_1}$ crosses strongly $X_{e_2}$ if and only if $G_{e_1}$ acts non-trivially on $T_2$
and $e_2\subset \Min_{T_2}(G_{e_1})$.
\end{lem}

\begin{proof}
Since the action of $G$ on $T_1$ induces a decomposition of $G$ as an amalgam or HNN extension,
$G$ is generated by elements sending $e_1$ to an edge having (at least) a common vertex with $e_1$. 
Thus, we can choose a finite generating set $S$ consisting of such elements.
For this choice of $S$, one gets 
$G_{e_1}=\partial X_{e_1}$:
any $g\in G_{e_1}$ lies in $\partial X_{e_1}$ because one can choose $s\in S\cup S\m$ sending $e_1$ into $\delta(e_1)^*$, 
which implies that $gs\notin X_{e_1}$;
for the other inclusion, if $g.e_1\in \delta(e_1)$ and $gs.e_1\in\delta(e_1)^*$, since $g.e_1$ and $gs.e_1$ are adjacent,
one necessarily gets $g.e_1=e_1$.

It follows that $X_{e_1}$ crosses strongly $X_{e_2}$ if and only if both $G_{e_1}\cap X_{e_2}$ and $G_{e_1}\cap X_{e_2}^*$ 
project to infinite sets in $C/G_{e_2}$. 
Since there is a natural bijection between $C/G_{e_2}$ and $G.e_2$,
this condition is equivalent to the fact that both $G_{e_1}.e_2\cap\delta(e_2)$ 
and $G_{e_1}.e_2\cap\delta(e_2)^*$ are infinite.

This clearly cannot occur if $G_{e_1}$ fixes a point $x$ in $T_2$: an isometry fixing $x$
cannot send $e_2$ to the component of $T\setminus \{e_2\}$ which does not contain $x$.
If $G_{e_1}$ does not fix a point in $T_2$ but every element is elliptic, then 
$G_{e_1}$ fixes an end $\omega$ of $T_2$, and a similar argument applies: 
an isometry fixing $\omega$ cannot send $e_2$ to the component of $T\setminus \{e_2\}$ which does not contain $\omega$ in its boundary.
Finally, if $G_{e_1}$ acts non-trivially on $T_2$, 
and if $e_2\notin \Min_{T_2}(G_{e_1})$, then similarly, no element of $G_{e_1}$ can send $e_2$ to
the component of  $T_2\setminus \{e_2\}$ which does not contain
$\Min_{T_2}(G_{e_1})$. Therefore, we proved that if $X_{e_1}$ crosses strongly $X_{e_2}$, then
$G_{e_1}$ acts non-trivially on $T_2$ and $e_2\in \Min_{T_2}(G_{e_1})$.

Conversely, assume that $e_2\in \Min_{T_2}(G_{e_1})$.
Then consider an element $h\in G_{e_1}$ whose axis in $T_2$ contains $e_2$.
Then clearly, $h^{\pm k}.e_2$ meets $\delta(e_2)$ and $\delta(e_2)^*$ infinitely many times.
\end{proof}

As a corollary, we immediately have the following interpretation of the strong intersection number.
This generalizes Corollary 3.16 in \cite{ScSw_splittings} which requires that the strong intersection number
coincide with the usual one.

\begin{SauveCompteurs}{cpt_carac_strong_intersection}
\begin{cor}\label{cor_carac_strong_intersection}
Consider $T_1,T_2$ two simplicial trees dual to one edge splittings a finitely generated group $G$.

Then  Scott and Swarup's (asymmetric) strong intersection number $si(T_1,T_2)$
is the number of two-cells of $\cala_1(T_1\times T_2)/G$.
This number can also be computed as the number of $G_{e_1}$-orbits of edges in
$\Min_{T_2}(G_{e_1})$ (where $e_1\in E(T_1)$ is any edge).
\end{cor}
\end{SauveCompteurs}

\section{Fujiwara and Papasoglu's enclosing groups}\label{sec_JSJ}

\begin{prop}[\cite{FuPa_JSJ}]\label{prop_JSJ_FP}
  Assume that $G$ acts minimally on two simplicial trees $T_1,T_2$, and that
for each $i\in\{1,2\}$,
\begin{itemize*}
\item  the stabilizer of each edge of $T_i$ acts non-trivially on the other tree;
\item  $G$ does not split over a subgroup of infinite index of the stabilizer of an edge of $T_i$.
\end{itemize*}

Then $\calc=\cala_1=\cala_2$.

Furthermore, assume that for each edge $e_i$ of $T_i$, the minimal $G_{e_i}$ invariant subtree of the other tree is a line.
Then $\calc \setminus V(\calc)$ is a surface. 
\end{prop}

\begin{proof}
Consider the asymmetric core $\cala_1=\cup_{x_1\in T_1}\{x_1\}\times \min_{T_2}(G_{x_1})$ defined above. 
We saw that $\cala_1\subset \calc$ and that $\cala_1$ is contractible.
The vertical fibers of $\cala_1$ are connected by definition.

We  prove that the horizontal fibers over an edge of $T_2$ is connected.
Assume that for some edge $e_2$, $F=p_2\m(e_2)\cap \cala_1$ is disconnected. First, there is at most one component $F$
which is $G_{e_2}$-invariant since such a component necessarily contains $\min_{T_1}(G_{e_2})\times \{e_2\}$. 
So consider another component $F_0=K_0\times\{e_2\}$ of $F$. The stabilizer $H$ of $F_0$ has infinite index in $G_{e_2}$
because $G_{e_2}$ acts non-trivially on $T_1$.
Since $\cala_1$ is simply connected, and since $F_0$ locally diconnects $\cala_1$ into two components, $F_0$ disconnects $\cala_1$ into two components.
This means that there is a tree $T_{F_0}$ dual to $G.F_0$ whose edges consist of connected components of $G.F_0$ and 
vertices to connected components of $\cala_1\setminus G.F_0$. This gives a splitting of $G$ over $H$. 
As soon as we prove that this splitting is non-trivial, we will get a contradiction showing that the horizontal fiber over $e_2$ is connected.

Since $G$ is finitely generated, we only need
to prove that for every vertex $v\in T_{F_0}$ there is an element $g\in G$ which does not fix $v$.
By equivariance, we can assume that $v$ is an endpoint of the edge $F_0$ of $T_{F_0}$.
Consider a point $x_1\in K_0$ and the corresponding an open edge $e=(x_1,e_2)$ in $F_0=K_o\times \{e_2\}$, and let
let $x_2\in T_2$ be the endpoint of $e_2$ such that $(x_1,x_2)$ is in the component $v$ of $\cala_1\setminus G.F_0$.
Let $g\in G_{x_1}$ whose axis contains $e_2$ (there exists such an element $g$ by definition of $\cala_1$).
Up to changing $g$ to its inverse, we can assume that $e_2\subset [x_2,g.x_2]$.
Now the path $\{x_1\}\times [x_2,g.x_2]$ in $\cala_1$ defines a path from $v$ to $g.v$ in $T_{F_0}$ which crosses exactly once
the edge $F_0$. This implies that $g.v\neq v$.

We now deduce that the horizontal fiber over a vertex $x_2\in T_2$ is connected.
Take two vertices $a,b\in p_2\m(x_2)$, and consider a piecewise linear path $c:[0,1]\ra \cala_1$ joining $a$ and $b$.
Let $c_2=p_2\circ c$. If $c_2$ is constant, we are done. Otherwise, consider a connected component $]s,t[$ of $[0,1]\setminus c_2\m({v_2})$.
Then $c_2(]s,s+\eps])$ and $c_2([t-\eps,t[)$ are contained in the same edge $e_2$.
Since $p_2\m(e_2)$ is connected, we can connect $c(s+\eps)$ to $c(t-\eps)$ by a path in $p_2\m(e_2)$,
and we can push it to $p_2\m(x_2)$ to get a path from $c(s)$ to $c(t)$ in $p_2\m(x_2)$.
Replacing $c_{|[s,t]}$ by such a path, and doing this for each connected component of $[0,1]\setminus c_2\m({x_2})$, we get
a path in $p_2\m(x_2)$ joining $a$ to $b$.

Since $\cala_1$ is connected and has connected fibers, we get $\cala_1\supset \calc$, and $\cala_1=\calc$.
The symmetric argument shows that $\cala_2=\calc$.

It is now easy to prove that $\calc\setminus V(\calc)$ is a surface under the additional hypothesis:
the link of each point in a horizontal edge $e_1 \times \{x_2\}$ in $\cala_1$ is a circle because
$\min_{T_2}(G_{x_1})$ is a line, and the symmetric argument shows that the link of a point in a vertical edge
is also a circle.
\end{proof}


\def\cprime{$'$}

\begin{flushleft}
Institut Fourier, UMR 5582,\\
BP74, Universit\'e Grenoble 1,\\
 38402 Saint-Martin d'H\`eres C\'edex, France.\\
\emph{e-mail}: \texttt{vincent.guirardel@ujf-grenoble.fr}
  
\end{flushleft}
\end{document}